# A Polynomial Chaos Expansion in Dependent Random Variables☆

Sharif Rahman[1]

*Applied Mathematical and Computational Sciences, The University of Iowa, Iowa City, Iowa 52242, U.S.A.*

**Abstract**

This paper introduces a new generalized polynomial chaos expansion (PCE) comprising measure-consistent multivariate orthonormal polynomials in dependent random variables. Unlike existing PCEs, whether classical or generalized, no tensor-product structure is assumed or required. Important mathematical properties of the generalized PCE are studied by constructing orthogonal decomposition of polynomial spaces, explaining completeness of orthogonal polynomials for prescribed assumptions, exploiting whitening transformation for generating orthonormal polynomial bases, and demonstrating mean-square convergence to the correct limit. Analytical formulae are proposed to calculate the mean and variance of a truncated generalized PCE for a general output variable in terms of the expansion coefficients. An example derived from a stochastic boundary-value problem illustrates the generalized PCE approximation in estimating the statistical properties of an output variable for 12 distinct non-product-type probability measures of input variables.

*Keywords:* Uncertainty quantification, multivariate orthogonal polynomials, Fourier series

## 1. Introduction

Polynomial chaos expansion (PCE) is an infinite series expansion of an output random variable involving orthogonal polynomials in input random variables. Introduced by Wiener [1] for Gaussian input variables, followed by a proof of convergence [2], the original PCE, referred to as the classical PCE in this paper, was later extended to a generalized PCE [3] to account for non-Gaussian variables. Approximations derived from truncated PCE, whether classical or generalized, are commonly used for solving uncertainty quantification problems, mostly in the context of solving stochastic partial differential equations [4, 5], yielding approximate second-moment statistics of an output random variable of interest. However, the existing PCE is largely founded on the independence assumption of input variables. The assumption exploits product-type probability measures, facilitating construction of the space of multivariate orthogonal polynomials via tensor product of the spaces of univariate orthogonal polynomials. In reality, there may exist significant correlation or dependence among input variables, impeding or invalidating many stochastic methods, including PCE. The Rosenblatt transformation [6], commonly used for mapping dependent to independent variables, may induce very strong nonlinearity to a stochastic response, potentially degrading or even prohibiting convergence of probabilistic solutions [7]. While the works of Soize and Ghanem [8] and Rahman [9] to cope with dependent variables are a step in the right direction, they, respectively, employ non-polynomial basis unamenable to producing analytical formulae for response statistics and focus strictly on Gaussian variables. Furthermore, the first of these studies does not address denseness or completeness of basis functions or account for infinitely many input variables. Therefore, innovations beyond tensor-product PCEs, capable of tackling non-product-type probability measures, are highly desirable.

☆Grant sponsor: U.S. National Science Foundation; Grant No. CMMI-1462385.
*Email address:* `sharif-rahman@uiowa.edu` (Sharif Rahman)
[1]Professor.

This study delves into a number of mathematical issues concerning necessary and sufficient conditions for the completeness of multivariate orthogonal polynomials; convergence, exactness, and optimal analyses; and approximation quality due to truncation – all associated with a generalized PCE for dependent, non-product-type probability measures. Therefore, the results of this paper are new in many aspects. The paper is organized as follows. Section 2 defines or discusses mathematical notations and preliminaries. A set of assumptions on the input probability measure required by the generalized PCE is explained. A brief exposition of multivariate orthogonal polynomials consistent with a general, non-product-type probability measure, including their second moment properties, is given in Section 3. The section also describes relevant polynomial spaces and construction of their orthogonal decompositions. The orthogonal basis and completeness of multivariate orthogonal polynomials have also been established. Section 4 defines the polynomial moment matrix, resulting in a variety of whitening transformations to produce measure-consistent orthonormal polynomials. The statistical properties of both orthogonal and orthonormal polynomials are presented. Section 5 formally introduces the generalized PCE for a square-integrable random variable. The convergence, exactness, and optimality of the generalized PCE are explained. In the same section, the approximation quality of a truncated generalized PCE is discussed. The formulae for the mean and variance of a truncated generalized PCE are derived, and methods for estimating the expansion coefficients are outlined. The section ends with an explanation on how and when the generalized PCE proposed can be extended for infinitely many input variables. The results from a simple yet illuminating example are reported in Section 6 with supplementary details in Appendix A. Finally, conclusions are drawn in Section 7.

## 2. Input Random Variables

Let $\mathbb{N} := \{1, 2, \ldots\}$, $\mathbb{N}_0 := \mathbb{N} \cup \{0\}$, $\mathbb{R} := (-\infty, +\infty)$, and $\mathbb{R}_0^+ := [0, +\infty)$ represent the sets of positive integer (natural), non-negative integer, real, and non-negative real numbers, respectively. For a non-zero, finite integer $N \in \mathbb{N}$, denote by $\mathbb{A}^N \subseteq \mathbb{R}^N$ a bounded or unbounded subdomain of $\mathbb{R}^N$. The set of $N \times N$ real-valued square matrices is denoted by $\mathbb{R}^{N\times N}$.

Let $(\Omega, \mathcal{F}, \mathbb{P})$ be a complete probability space, where $\Omega$ is a sample space representing an abstract set of elementary events, $\mathcal{F}$ is a $\sigma$-algebra on $\Omega$, and $\mathbb{P} : \mathcal{F} \to [0, 1]$ is a probability measure. With $\mathcal{B}^N := \mathcal{B}(\mathbb{A}^N)$ representing the Borel $\sigma$-algebra on $\mathbb{A}^N \subseteq \mathbb{R}^N$, consider an $\mathbb{A}^N$-valued input random vector $\mathbf{X} := (X_1, \ldots, X_N)^T : (\Omega, \mathcal{F}) \to (\mathbb{A}^N, \mathcal{B}^N)$, describing the statistical uncertainties in all system parameters of a stochastic problem. The input random variables are also referred to as basic random variables. The integer $N$ represents the number of input random variables and is referred to as the dimension of the stochastic problem.

Denote by $F_{\mathbf{X}}(\mathbf{x}) := \mathbb{P}(\cap_{i=1}^N \{X_i \le x_i\})$ the joint distribution function of $\mathbf{X}$, admitting the joint probability density function $f_{\mathbf{X}}(\mathbf{x}) := \partial^N F_{\mathbf{X}}(\mathbf{x})/\partial x_1 \cdots \partial x_N$. Given the abstract probability space $(\Omega, \mathcal{F}, \mathbb{P})$, the image probability space is $(\mathbb{A}^N, \mathcal{B}^N, f_{\mathbf{X}} d\mathbf{x})$, where $\mathbb{A}^N$ can be viewed as the image of $\Omega$ from the mapping $\mathbf{X} : \Omega \to \mathbb{A}^N$, and is also the support of $f_{\mathbf{X}}(\mathbf{x})$. Relevant statements and objects in one space have obvious counterparts in the other space. Both probability spaces will be used in this paper.

A set of assumptions used for or required by the generalized PCE is as follows.

**Assumption 1.** *The random vector* $\mathbf{X} := (X_1, \ldots, X_N)^T : (\Omega, \mathcal{F}) \to (\mathbb{A}^N, \mathcal{B}^N)$

(1) *has an absolutely continuous joint distribution function* $F_{\mathbf{X}}(\mathbf{x})$ *and a continuous joint probability density function* $f_{\mathbf{X}}(\mathbf{x})$ *with a bounded or unbounded support* $\mathbb{A}^N \subseteq \mathbb{R}^N$*;*

(2) *possesses absolute finite moments of all orders, that is, for all* $\mathbf{j} := (j_1, \ldots, j_N) \in \mathbb{N}_0^N$*,*

$$\mathbb{E}\left[|\mathbf{X}^{\mathbf{j}}|\right] := \int_{\Omega} |\mathbf{X}(\omega)|^{\mathbf{j}} d\mathbb{P}(\omega) = \int_{\mathbb{A}^N} |\mathbf{x}^{\mathbf{j}}| f_{\mathbf{X}}(\mathbf{x}) d\mathbf{x} < \infty,$$

*where* $\mathbf{X}^{\mathbf{j}} = X_1^{j_1} \cdots X_N^{j_N}$ *and* $\mathbb{E}$ *is the expectation operator with respect to the probability measure* $\mathbb{P}$ *or* $f_{\mathbf{X}}(\mathbf{x})d\mathbf{x}$*; and*

(3) *has a joint probability density function* $f_{\mathbf{X}}(\mathbf{x})$*, which*

*(a) has a compact support, that is, there exists a compact subset $\mathbb{A}^N \subset \mathbb{R}^N$ such that $\mathbb{P}(\mathbf{X} \in \mathbb{A}^N) = 1$, or*

*(b) is exponentially integrable, that is, there exists a real number $\alpha > 0$ such that*

$$\int_{\mathbb{A}^N} \exp\left(\alpha \|\mathbf{x}\|\right) f_{\mathbf{X}}(\mathbf{x}) d\mathbf{x} < \infty,$$

*where $\|\cdot\| : \mathbb{A}^N \to \mathbb{R}_0^+$ is an arbitrary norm.*

Item (1) of Assumption 1 is not essential to PCE, but it is commonly invoked in applications. Item (2) of Assumption 1 assures the existence of an infinite sequence of multivariate orthogonal polynomials consistent with the input probability measure. Item (3) of Assumption 1, in addition to Items (1) and (2), guarantees the input probability measure to be determinate[2], resulting in a complete orthogonal polynomial system and hence a basis of a function space of interest. The assumptions impose only mild restrictions on the probability measure. Examples of input random variables satisfying Assumption 1 are multivariate Gaussian, uniform, exponential, Laplace variables, including some endowed with rotationally invariant density functions [10]. This assumption, to be explained in the next section, is vitally important for the determinacy of the probability measure and the completeness of orthogonal polynomials. Examples where Items (1) and (2) are satisfied, but Item (3) is not, are lognormal distributions, select distributions from the Farlie-Gumbel-Morgenstern family, Kotz-type distributions, and cases involving nonlinear transformations of random variables with determinate distributions [11]. As noted by Ernst et al. [12] for the univariate case ($N = 1$), the violation of Item (3) leads to indeterminacy of the lognormal probability measure and thereby fails to form a complete orthogonal polynomial system.

## 3. Multivariate Orthogonal Polynomials

Let $\mathbf{j} := (j_1, \ldots, j_N) \in \mathbb{N}_0^N$, $j_i \in \mathbb{N}_0$, denote an $N$-dimensional multi-index. For $\mathbf{x} = (x_1, \ldots, x_N) \in \mathbb{A}^N \subseteq \mathbb{R}^N$, a monomial in the variables $x_1, \ldots, x_N$ is the product $\mathbf{x}^{\mathbf{j}} = x_1^{j_1} \cdots x_N^{j_N}$ and has a total degree $|\mathbf{j}| = j_1 + \cdots + j_N$. A linear combination of $\mathbf{x}^{\mathbf{j}}$, where $|\mathbf{j}| = l$, $l \in \mathbb{N}_0$, is a homogeneous polynomial of degree $l$. Denote by

$$\mathcal{P}_l^N := \operatorname{span}\{\mathbf{x}^{\mathbf{j}} : |\mathbf{j}| = l, \mathbf{j} \in \mathbb{N}_0^N\}, \ l \in \mathbb{N}_0,$$

the space of homogeneous polynomials of degree $l$, by

$$\Pi_m^N := \operatorname{span}\{\mathbf{x}^{\mathbf{j}} : |\mathbf{j}| \le m, \mathbf{j} \in \mathbb{N}_0^N\}, \ m \in \mathbb{N}_0,$$

the space of polynomials of degree at most $m$, and by

$$\Pi^N := \mathbb{R}[\mathbf{x}] = \mathbb{R}[x_1, \ldots, x_N],$$

the space of all real polynomials in $\mathbf{x}$. It is well known that the dimensions of the vector spaces $\mathcal{P}_l^N$ and $\Pi_m^N$, respectively, are [10]

$$\dim \mathcal{P}_l^N = \#\left\{\mathbf{j} \in \mathbb{N}_0^N : |\mathbf{j}| = l\right\} = \binom{N+l-1}{l} =: K_{N,l} \tag{1}$$

and

$$\dim \Pi_m^N = \sum_{l=0}^{m} \dim \mathcal{P}_l^N = \sum_{l=0}^{m} \binom{N+l-1}{l} = \binom{N+m}{m}.$$

[2]The density function of the probability measure, if it is uniquely determined by a sequence of moments, is called determinate or M-determinate. Otherwise, the density function is indeterminate or M-indeterminate. This is known as the moment problem with three prominent flavors, depending on the support of the density: Hausdorff moment problem ($\mathbb{A}^N = [0,1]^N$), Stieltjes moment problem ($\mathbb{A}^N = \mathbb{R}_0^{+N}$), and Hamburger moment problem ($\mathbb{A}^N = \mathbb{R}^N$).

*3.1. Orthogonal polynomials*

Let $f_{\mathbf{X}}(\mathbf{x})d\mathbf{x}$ be a probability measure on $\mathbb{A}^N$, satisfying Assumption 1. For any polynomial pair $P, Q \in \Pi^N$, define an inner product

$$(P,Q)_{f_{\mathbf{X}}d\mathbf{x}} := \int_{\mathbb{A}^N} P(\mathbf{x})Q(\mathbf{x})f_{\mathbf{X}}(\mathbf{x})d\mathbf{x} = \mathbb{E}\left[P(\mathbf{X})Q(\mathbf{X})\right] \tag{2}$$

on $\Pi^N$ with respect to the measure $f_{\mathbf{X}}(\mathbf{x})d\mathbf{x}$ and the induced norm

$$\|P\|_{f_{\mathbf{X}}d\mathbf{x}} := \sqrt{(P,P)_{f_{\mathbf{X}}d\mathbf{x}}} = \left(\int_{\mathbb{A}^N} P^2(\mathbf{x})f_{\mathbf{X}}(\mathbf{x})d\mathbf{x}\right)^{1/2} = \sqrt{\mathbb{E}\left[P^2(\mathbf{X})\right]}.$$

The polynomials $P \in \Pi^N$ and $Q \in \Pi^N$ are called orthogonal to each other with respect to $f_{\mathbf{X}}(\mathbf{x})d\mathbf{x}$ if $(P,Q)_{f_{\mathbf{X}}d\mathbf{x}} = 0$. Moreover, a polynomial $P \in \Pi^N$ is said to be an orthogonal polynomial with respect to $f_{\mathbf{X}}(\mathbf{x})d\mathbf{x}$ if it is orthogonal to all polynomials of lower degree, that is, if [10]

$$(P,Q)_{f_{\mathbf{X}}d\mathbf{x}} = 0 \ \forall Q \in \Pi^N \text{ with } \deg Q < \deg P. \tag{3}$$

Under Items (1) and (2) of Assumption 1, moments of $\mathbf{X}$ of all orders exist and are finite, so that the inner product in (2) is well defined. As the inner product is positive-definite, clearly $\|P\|_{f_{\mathbf{X}}d\mathbf{x}} > 0$ for all non-zero $P \in \Pi^N$. Then there exists an infinite set of multivariate orthogonal polynomials, say, $\{P_{\mathbf{j}}(\mathbf{x}) : \mathbf{j} \in \mathbb{N}_0^N\}$, $P_{\mathbf{0}} = 1$, $P_{\mathbf{j}} \neq 0$, which is consistent with the probability measure $f_{\mathbf{X}}(\mathbf{x})d\mathbf{x}$, satisfying

$$(P_{\mathbf{j}}, P_{\mathbf{k}})_{f_{\mathbf{X}}d\mathbf{x}} = 0 \text{ whenever } |\mathbf{j}| \neq |\mathbf{k}| \tag{4}$$

for $\mathbf{k} \in \mathbb{N}_0^N$. Here, the multi-index $\mathbf{j}$ of the multivariate polynomial $P_{\mathbf{j}}(\mathbf{x})$ refers to its total degree $|\mathbf{j}| = j_1 + \cdots + j_N$. Clearly, each $P_{\mathbf{j}} \in \Pi^N$ is an orthogonal polynomial satisfying (3). This means that $P_{\mathbf{j}}$ is orthogonal to all polynomials of different degrees, but it may not be orthogonal to other orthogonal polynomials of the same degree.

Consider for each $l \in \mathbb{N}_0$ the elements of the set $\{\mathbf{j} \in \mathbb{N}_0^N : |\mathbf{j} = l\}$, $l \in \mathbb{N}_0$, which is arranged as $\mathbf{j}^{(1)}, \ldots, \mathbf{j}^{(K_{N,l})}$ according to a monomial order of choice. The set has cardinality $K_{N,l}$ as defined in (1). Denote by

$$\mathbf{x}_l = (\mathbf{x}^{\mathbf{j}^{(1)}}, \ldots, \mathbf{x}^{\mathbf{j}^{(K_{N,l})}})^T$$

the $K_{N,l}$-dimensional column vector whose elements are the monomials $\mathbf{x}^{\mathbf{j}}$ for $|\mathbf{j}| = l$ and by

$$\mathbf{P}_l(\mathbf{x}) := (P_{l,\mathbf{j}^{(1)}}(\mathbf{x}), \ldots, P_{l,\mathbf{j}^{(K_{N,l})}}(\mathbf{x}))^T \tag{5}$$

the $K_{N,l}$-dimensional column vector whose elements are obtained from the polynomial sequence $\{P_{l,\mathbf{j}}(\mathbf{x}) := P_{\mathbf{j}}(\mathbf{x})\}_{|\mathbf{j}|=l}$, both arranged in the aforementioned order. This leads to a formal definition of multivariate orthogonal polynomials.

**Definition 2** (Dunkl and Xu [10])**.** *Let $(\cdot,\cdot)_{f_{\mathbf{X}}d\mathbf{x}} : \Pi^N \times \Pi^N \to \mathbb{R}$ be an inner product. A set of polynomials $\{P_{\mathbf{j}}(\mathbf{x}) : |\mathbf{j}| = l, \mathbf{j} \in \mathbb{N}_0^N\}$, $P_{\mathbf{j}}(\mathbf{x}) \in \Pi_l^N$, of degree $l$ or its $K_{N,l}$-dimensional column vector $\mathbf{P}_l(\mathbf{x})$, is said to be orthogonal with respect to the inner product $(\cdot,\cdot)_{f_{\mathbf{X}}d\mathbf{x}}$, or alternatively with respect to the probability measure $f_{\mathbf{X}}(\mathbf{x})d\mathbf{x}$, if, for $l, r \in \mathbb{N}_0$,*

$$\left(\mathbf{x}_r, \mathbf{P}_l^T(\mathbf{x})\right)_{f_{\mathbf{X}}d\mathbf{x}} := \int_{\mathbb{A}^N} \mathbf{x}_r \mathbf{P}_l^T(\mathbf{x}) f_{\mathbf{X}}(\mathbf{x})d\mathbf{x} =: \mathbb{E}\left[\mathbf{X}_r \mathbf{P}_l^T(\mathbf{X})\right] = \mathbf{0}, \quad l > r, \tag{6}$$

*where*

$$\mathbf{S}_l := \left(\mathbf{x}_l, \mathbf{P}_l^T(\mathbf{x})\right)_{f_{\mathbf{X}}d\mathbf{x}} := \int_{\mathbb{A}^N} \mathbf{x}_l \mathbf{P}_l^T(\mathbf{x}) f_{\mathbf{X}}(\mathbf{x})d\mathbf{x} =: \mathbb{E}\left[\mathbf{X}_l \mathbf{P}_l^T(\mathbf{X})\right] \tag{7}$$

*is a $K_{N,l} \times K_{N,l}$ invertible matrix.*

Using the vector notation, one can write

$$\mathbf{P}_r(\mathbf{x}) = \mathbf{H}_{r,r}\mathbf{x}_r + \mathbf{H}_{r,r-1}\mathbf{x}_{r-1} + \cdots + \mathbf{H}_{r,0}\mathbf{x}_0, \quad r \in \mathbb{N}_0,$$

where $\mathbf{H}_{r,r-k}$, $k = 0, 1, \ldots, r$, are various coefficient matrices of size $K_{N,r} \times K_{N,r-k}$. Then, using (6) and (7), the inner products $(\mathbf{P}_r(\mathbf{x}), \mathbf{P}_l^T(\mathbf{x}))_{f_{\mathbf{X}}d\mathbf{x}} = \mathbf{0}$ when $l > r$ and $(\mathbf{P}_r(\mathbf{x}), \mathbf{P}_l^T(\mathbf{x}))_{f_{\mathbf{X}}d\mathbf{x}} = \mathbf{H}_{l,l}\mathbf{S}_l$ when $l = r$. Therefore, Definition 2 agrees with the usual notion of orthogonal polynomials satisfying (4). Perhaps the most prominent example of classical multivariate orthogonal polynomials is the case of multivariate Hermite polynomials, which are consistent with the measure defined by a Gaussian density on $\mathbb{R}^N$ [13, 14]. Readers interested to learn more about orthogonal polynomials in multiple variables are referred to the works of Appell and de Fériet [15], Erdelyi [13], Krall and Sheffer [16], and Dunkl and Xu [10].

For general probability measures, established numerical techniques, such as the Gram-Schmidt orthogonalization process [17], can be applied to a sequence of monomials $\{\mathbf{x}^{\mathbf{j}}\}_{\mathbf{j}\in\mathbb{N}_0^N}$ with respect to the inner product in (2) to generate a corresponding sequence of any measure-consistent orthogonal polynomials. However, an important difference between univariate polynomials and multivariate polynomials is the lack of an obvious natural order in the latter. The natural order for monomials of univariate polynomials is the degree order; that is, one orders monomials according to their degree. For multivariate polynomials, there are many options, such as lexicographic order, graded lexicographic order, and reversed graded lexicographic order, to name just three. There is no natural choice, and different orders will give different sequences of orthogonal polynomials from the Gram-Schmidt process. It is important to emphasize that the space of multivariate orthogonal polynomials for a generally non-product-type density function cannot be constructed by the tensor product of the spaces of univariate orthogonal polynomials.

### *3.2. Orthogonal decomposition of polynomial spaces*

Let $\mathcal{V}_0^N := \Pi_0^N = \text{span}\{1\}$ be the space of constant functions. For each $1 \le l < \infty$, denote by $\mathcal{V}_l^N \subset \Pi_l^N$ the space of orthogonal polynomials of degree exactly $l$ that are orthogonal to all polynomials in $\Pi_{l-1}^N$, that is,

$$\mathcal{V}_l^N := \{P \in \Pi_l^N : (P, Q)_{f_{\mathbf{X}}d\mathbf{x}} = 0 \; \forall \, Q \in \Pi_{l-1}^N\}, \; 1 \le l < \infty.$$

Then $\mathcal{V}_l^N$, provided that the support of $f_{\mathbf{X}}(\mathbf{x})$ has a non-empty interior, is a vector space of dimension [10]

$$K_{N,l} := \dim \mathcal{V}_l^N = \dim \mathcal{P}_l^N = \binom{N+l-1}{l}.$$

Many choices exist for the basis of $\mathcal{V}_l^N$; the bases of $\mathcal{V}_l^N$ do not have to be mutually orthogonal. Furthermore, with the exception of the monic orthogonal polynomials, the bases are not unique in the multivariate case. Here, to be formally proved in the next section, select $\{P_{\mathbf{j}}(\mathbf{x}) : |\mathbf{j}| = l, \mathbf{j} \in \mathbb{N}_0^N\} \subset \mathcal{V}_l^N$ to be a basis of $\mathcal{V}_l^N$, comprising $K_{N,l}$ number of basis functions. Each basis function $P_{\mathbf{j}}(\mathbf{x})$ is a multivariate orthogonal polynomial of degree $|\mathbf{j}|$ as discussed earlier. Obviously,

$$\mathcal{V}_l^N = \text{span}\{P_{\mathbf{j}} : |\mathbf{j}| = l, \mathbf{j} \in \mathbb{N}_0^N\}, \; 0 \le l < \infty.$$

According to (4), $P_{\mathbf{j}}$ is orthogonal to $P_{\mathbf{k}}$ whenever $|\mathbf{j}| \neq |\mathbf{k}|$. Therefore, any two polynomial subspaces $\mathcal{V}_l^N$ and $\mathcal{V}_r^N$, where $0 \le l, r < \infty$, are orthogonal whenever $l \neq r$. In consequence, there exist orthogonal decompositions of

$$\Pi_m^N = \bigoplus_{l=0}^{m} \mathcal{V}_l^N = \bigoplus_{l=0}^{m} \text{span}\{P_{\mathbf{j}} : |\mathbf{j}| = l, \mathbf{j} \in \mathbb{N}_0^N\} = \text{span}\{P_{\mathbf{j}} : 0 \le |\mathbf{j}| \le m, \mathbf{j} \in \mathbb{N}_0^N\}$$

and

$$\Pi^N = \bigoplus_{l\in\mathbb{N}_0} \mathcal{V}_l^N = \bigoplus_{l\in\mathbb{N}_0} \text{span}\{P_{\mathbf{j}} : |\mathbf{j}| = l, \mathbf{j} \in \mathbb{N}_0^N\} = \text{span}\{P_{\mathbf{j}} : \mathbf{j} \in \mathbb{N}_0^N\} \tag{8}$$

with the symbol $\oplus$ representing orthogonal sum.

*3.3. Completeness of orthogonal polynomials and basis*

An important question regarding orthogonal polynomials is whether they are complete and constitute a basis in a function space of interest, such as a Hilbert space. Let $L^2(\mathbb{A}^N, \mathcal{B}^N, f_{\mathbf{X}}d\mathbf{x})$ represent a Hilbert space of square-integrable functions with respect to the probability measure $f_{\mathbf{X}}(\mathbf{x})d\mathbf{x}$ supported on $\mathbb{A}^N$. The following two propositions show that, indeed, orthogonal polynomials span various spaces of interest.

**Proposition 3.** *Let $\mathbf{X} := (X_1, \ldots, X_N)^T : (\Omega, \mathcal{F}) \to (\mathbb{A}^N, \mathcal{B}^N)$, $N \in \mathbb{N}$, be an $N$-dimensional random vector with multivariate probability density function $f_{\mathbf{X}}(\mathbf{x})$, satisfying Assumption 1. Then $\{P_{\mathbf{j}}(\mathbf{x}) : |\mathbf{j}| = l, \mathbf{j} \in \mathbb{N}_0^N\}$, the set of multivariate orthogonal polynomials of degree $l$ consistent with the probability measure $f_{\mathbf{X}}(\mathbf{x})d\mathbf{x}$, is a basis of $\mathcal{V}_l^N$.*

*Proof.* Under Items (1) and (2) of Assumption 1, orthogonal polynomials with respect to the probability measure $f_{\mathbf{X}}(\mathbf{x})d\mathbf{x}$ exist. Let $\mathbf{a}_l^T = (a_{l,1}, \ldots, a_{l,K_{N,l}})$ be a row vector comprising some constants $a_{l,i} \in \mathbb{R}$, $i = 1, \ldots, K_{N,l}$. Set $\mathbf{a}_l^T \mathbf{P}_l(\mathbf{x}) = 0$. Multiply both sides of the equality from the right by $\mathbf{x}_l^T$, integrate with respect to the measure $f_{\mathbf{X}}(\mathbf{x})d\mathbf{x}$ over $\mathbb{A}^N$, and apply transposition to obtain

$$\mathbf{S}_l \mathbf{a}_l = \mathbf{0}, \tag{9}$$

where $\mathbf{S}_l$, defined in (7), is a $K_{N,l} \times K_{N,l}$ invertible matrix. Therefore, (9) yields $\mathbf{a}_l = \mathbf{0}$, proving linear independence of the elements of $\mathbf{P}_l(\mathbf{x})$ or the set $\{P_{\mathbf{j}}(\mathbf{x}) : |\mathbf{j}| = l, \mathbf{j} \in \mathbb{N}_0^N\}$. Furthermore, the dimension $K_{N,l}$ of $\mathcal{V}_l^N$ matches exactly the number of elements of the aforementioned set. Therefore, the spanning set $\{P_{\mathbf{j}}(\mathbf{x}) : |\mathbf{j}| = l, \mathbf{j} \in \mathbb{N}_0^N\}$ forms a basis of $\mathcal{V}_l^N$. □

**Proposition 4.** *Let $\mathbf{X} := (X_1, \ldots, X_N)^T : (\Omega, \mathcal{F}) \to (\mathbb{A}^N, \mathcal{B}^N)$, $N \in \mathbb{N}$, be an $N$-dimensional random vector with multivariate probability density function $f_{\mathbf{X}}(\mathbf{x})$, satisfying Assumption 1. Consistent with the probability measure $f_{\mathbf{X}}(\mathbf{x})d\mathbf{x}$, let $\{P_{\mathbf{j}}(\mathbf{x}) : |\mathbf{j}| = l, \mathbf{j} \in \mathbb{N}_0^N\}$, the set of multivariate orthogonal polynomials of degree $l$, be a basis of $\mathcal{V}_l^N$. Then the set of polynomials from the orthogonal sum*

$$\bigoplus_{l \in \mathbb{N}_0} \mathit{span}\{P_{\mathbf{j}}(\mathbf{x}) : |\mathbf{j}| = l, \mathbf{j} \in \mathbb{N}_0^N\}$$

*is dense in $L^2(\mathbb{A}^N, \mathcal{B}^N, f_{\mathbf{X}}d\mathbf{x})$. Moreover,*

$$L^2(\mathbb{A}^N, \mathcal{B}^N, f_{\mathbf{X}}d\mathbf{x}) = \overline{\bigoplus_{l \in \mathbb{N}_0} \mathcal{V}_l^N} \tag{10}$$

*where the overline denotes set closure.*

*Proof.* Under Items (1) and (2) of Assumption 1, orthogonal polynomials with respect to the probability measure $f_{\mathbf{X}}(\mathbf{x})d\mathbf{x}$ exist. According to Theorem 3.2.18 of Dunkl and Xu [10] and related discussion, which exploits Items 3(a) and 3(b) of Assumption 1, the polynomial space $\Pi^N$ is dense in the space $L^2(\mathbb{A}^N, \mathcal{B}^N, f_{\mathbf{X}}d\mathbf{x})$. Therefore, the set of polynomials from the orthogonal sum, which is equal to $\Pi^N$ as per (8), is dense in $L^2(\mathbb{A}^N, \mathcal{B}^N, f_{\mathbf{X}}d\mathbf{x})$. Including the limit points of the orthogonal sum yields (10). □

## 4. Multivariate orthonormal polynomials

Once the multivariate orthogonal polynomials are obtained, they can be linearly transformed to generate multivariate orthonormal polynomials. The latter polynomials, while they are not required [9], result in concise forms of the generalized PCE and second-moment properties of an output random variable of interest.

*4.1. Polynomial moment matrix*

When the input random variables $X_1, \ldots, X_N$, instead of the variables $x_1, \ldots, x_N$, are inserted in the argument, $\mathbf{P}_l$ in (5) becomes a vector of random orthogonal polynomials. A formal definition of the polynomial moment matrix follows.

**Definition 5.** *Let* $\mathbf{P}_l(\mathbf{X}) := (P_{l,\mathbf{j}^{(1)}}(\mathbf{X}), \ldots, P_{l,\mathbf{j}^{(K_{N,l})}}(\mathbf{X}))^T$*,* $l \in \mathbb{N}_0$*, be a* $K_{N,l}$*-dimensional vector of constant or random orthogonal polynomials. The* $K_{N,l} \times K_{N,l}$ *matrix, defined by*

$$\mathbf{G}_l := \mathbb{E}[\mathbf{P}_l(\mathbf{X})\mathbf{P}_l^T(\mathbf{X})], \tag{11}$$

*with its* $(p,q)$*th element*

$$G_{l,pq} = \mathbb{E}[P_{l,\mathbf{j}^{(p)}}(\mathbf{X})P_{l,\mathbf{j}^{(q)}}(\mathbf{X})], \ p,q = 1,\ldots,K_{N,l},$$

*is called the polynomial moment matrix of* $\mathbf{P}_l(\mathbf{X})$*.*

When $l=0$, $K_{N,0}=1$ and $\mathbf{P}_0(\mathbf{X}) = (1)^T = 1$ regardless of $N$. Therefore, (11) from Definition 5 yields $\mathbf{G}_0 = [1]$ to be a $1 \times 1$ matrix. When $l > 0$, $G_{l,pq}$ represents the covariance between two random polynomials of degree $l$, as it will be shown later that $\mathbb{E}[\mathbf{P}_l(\mathbf{X})] = \mathbf{0}$ when $l > 0$. In this case, $\mathbf{G}_l$ is nothing but the covariance matrix of $\mathbf{P}_l(\mathbf{X})$.

**Proposition 6.** *The polynomial moment matrix* $\mathbf{G}_l$ *is symmetric and positive-definite.*

*Proof.* By definition, $\mathbf{G}_l = \mathbf{G}_l^T$. From Proposition 3, the elements of $\mathbf{P}_l(\mathbf{x})$ are linearly independent. Therefore, for any $\mathbf{0} \neq \boldsymbol{\alpha}_l \in \mathbb{R}^{K_{N,l}}$, $\boldsymbol{\alpha}_l^T\mathbf{P}_l(\mathbf{x}) \in \Pi^N$ is a non-zero polynomial, satisfying

$$\boldsymbol{\alpha}_l^T\mathbf{G}_l\boldsymbol{\alpha}_l = \mathbb{E}\left[\left(\boldsymbol{\alpha}_l^T\mathbf{P}_l(\mathbf{X})\right)^2\right] = \|\boldsymbol{\alpha}_l^T\mathbf{P}_l(\mathbf{x})\|^2_{f_{\mathbf{X}}d\mathbf{x}} > 0,$$

as the inner product defined in (2) is positive-definite on $\Pi^N$. Therefore, $\mathbf{G}_l$ is a symmetric, positive-definite matrix. □

*4.2. Whitening transformation*

From Proposition 6, $\mathbf{G}_l$ is positive-definite and hence invertible. Therefore, for each $l \in \mathbb{N}_0$, there exists a non-singular matrix $\mathbf{W}_l \in \mathbb{R}^{K_{N,l} \times K_{N,l}}$ such that

$$\mathbf{W}_l^T\mathbf{W}_l = \mathbf{G}_l^{-1} \ \text{ or } \ \mathbf{W}_l^{-1}\mathbf{W}_l^{-T} = \mathbf{G}_l. \tag{12}$$

This leads to multivariate orthonormal polynomials as follows.

**Definition 7.** *Let* $\mathbf{X} := (X_1,\ldots,X_N)^T$ *be a vector of* $N \in \mathbb{N}$ *input random variables fulfilling Assumption 1. Then, given the vector* $\mathbf{P}_l(\mathbf{x}) \in \mathbb{R}^{K_{N,l}}$ *of multivariate orthogonal polynomials of degree* $l$*, the corresponding vector*

$$\boldsymbol{\Psi}_l(\mathbf{x}) := (\Psi_{l,\mathbf{j}^{(1)}}(\mathbf{x}), \ldots, \Psi_{l,\mathbf{j}^{(K_{N,l})}}(\mathbf{x}))^T \in \mathbb{R}^{K_{N,l}}$$

*of multivariate orthonormal polynomials, also of degree* $l$*, is obtained from the whitening transformation*

$$\boldsymbol{\Psi}_l(\mathbf{x}) = \mathbf{W}_l\mathbf{P}_l(\mathbf{x}), \ l \in \mathbb{N}_0, \tag{13}$$

*where* $\mathbf{W}_l \in \mathbb{R}^{K_{N,l} \times K_{N,l}}$ *is a non-singular whitening matrix satisfying* (12)*.*

The whitening transformation in Definition 7 is a linear transformation that converts $\mathbf{P}_l(\mathbf{x})$ into $\boldsymbol{\Psi}_l(\mathbf{x})$ in such a way that the latter has uncorrelated random polynomials, each with *zero* means for $l > 0$. The transformation is called "whitening" because it changes one random vector to the other, which has statistical properties akin to that of a white noise vector. However, the condition (12) does not uniquely determine the whitening matrix $\mathbf{W}_l$. Indeed, there exists infinitely many choices of $\mathbf{W}_l$ that all satisfy (12). All of these choices result in a linear transformation, decorrelating $\mathbf{P}_l(\mathbf{x})$ but producing different random vectors $\boldsymbol{\Psi}_l(\mathbf{x})$. As demonstrated by Kessy et al. [18], the selection of $\mathbf{W}_l$ depends on the desired cross-covariance or cross-correlation between $\mathbf{P}_l(\mathbf{x})$ and $\boldsymbol{\Psi}_l(\mathbf{x})$. Table 1 lists five commonly used whitening matrices from practical applications. Any of these whitening matrices and possibly others, in conjunction with (13), can be used to generate multivariate orthonormal polynomials.

The whitening transformation should not be confused with measure transformations commonly used for decorrelating dependent Gaussian variables with positive-definite covariance matrices. Such transformations are generally nonlinear for non-Gaussian probability measures. In contrast, the transformation introduced here is linear and decorrelates instead random orthogonal polynomials for any probability measure of $\mathbf{X}$. Therefore, a wide variety of input variables, including non-Gaussian variables, can be dealt with when generating measure-consistent orthonormal polynomials.

Table 1: Five choices for the whitening matrix $\mathbf{W}_l$.

| Name | Whitening matrix $\mathbf{W}_l$ | Notes |
|---|---|---|
| Zero − phase component analysis (ZCA) | $\mathbf{G}_l^{-1/2}$ | $\mathbf{G}_l^{-1/2} = \mathbf{U}_l \boldsymbol{\Lambda}_l^{-1/2} \mathbf{U}_l^T$, where $\mathbf{U}_l$ and $\boldsymbol{\Lambda}_l$ contain eigenvectors and eigenvalues of $\mathbf{G}_l$, respectively. |
| Principal component analysis (PCA) | $\boldsymbol{\Lambda}_l^{-1/2}\mathbf{U}_l^T$ | See notes in the second row from the top. |
| Cholesky decomposition | $\mathbf{L}_l^T$ | $\mathbf{G}_l^{-1} = \mathbf{L}_l \mathbf{L}_l^T$. |
| ZCA − correlation adjusted | $\bar{\mathbf{G}}_l^{-1/2}\mathbf{V}_l^{-1/2}$ | $\mathbf{G}_l = \mathbf{V}_l^{1/2} \bar{\mathbf{G}}_l \mathbf{V}_l^{1/2}$. |
| PCA − correlation adjusted | $\bar{\boldsymbol{\Lambda}}_l^{-1/2}\bar{\mathbf{U}}_l^T\mathbf{V}_l^{-1/2}$ | $\bar{\mathbf{G}}_l^{-1/2} = \bar{\mathbf{U}}_l \bar{\boldsymbol{\Lambda}}_l^{-1/2} \bar{\mathbf{U}}_l^T$, where $\bar{\mathbf{U}}_l$ and $\bar{\boldsymbol{\Lambda}}_l$ contain eigenvectors and eigenvalues of $\bar{\mathbf{G}}_l$, respectively; $\mathbf{G}_l = \mathbf{V}_l^{1/2} \bar{\mathbf{G}}_l \mathbf{V}_l^{1/2}$. |

### *4.3. Statistical properties*

The polynomial vectors $\mathbf{P}_l(\mathbf{X})$ and $\boldsymbol{\Psi}_l(\mathbf{X})$ are functions of random input variables. Therefore, it is important to establish their second-moment properties, to be exploited in Sections 5 and 6.

**Proposition 8.** *Let $\mathbf{X} := (X_1, \dots, X_N)^T$ be a vector of $N \in \mathbb{N}$ input random variables fulfilling Assumption 1. For $l, r \in \mathbb{N}_0$, the first- and second-order moments of the vector of multivariate orthogonal polynomials are*

$$\mathbb{E}\left[\mathbf{P}_l(\mathbf{X})\right] = \begin{cases} (1)^T = (1), & l = 0, \\ \mathbf{0}, & l \neq 0. \end{cases} \tag{14}$$

*and*

$$\mathbb{E}\left[\mathbf{P}_l(\mathbf{X})\mathbf{P}_r^T(\mathbf{X})\right] = \begin{cases} \mathbf{G}_l, & l = r, \\ \mathbf{0}, & l \neq r. \end{cases} \tag{15}$$

*respectively, where $\mathbf{G}_l \in \mathbb{R}^{K_{N,l} \times K_{N,l}}$ is defined by* (11).

*Proof.* The non-trivial result of (14) is attained from the recognition that $P_{\mathbf{0}} = 1$, whereas the trivial result of (14) follows by setting $|\mathbf{j}| = l \neq 0$ and $|\mathbf{k}| = 0$ in (4). The non-trivial and trivial results of (15) are obtained, respectively, by (11) of Definition 5 and by setting $|\mathbf{j}| = l \neq 0$ and $|\mathbf{k}| = r$, $l \neq r$, in (4). □

**Proposition 9.** *Let $\mathbf{X} := (X_1, \dots, X_N)^T$ be a vector of $N \in \mathbb{N}$ input random variables fulfilling Assumption 1. For $l, r \in \mathbb{N}_0$, the first- and second-order moments of the vector of multivariate orthonormal polynomials are*

$$\mathbb{E}\left[\boldsymbol{\Psi}_l(\mathbf{X})\right] = \begin{cases} (1)^T = (1), & l = 0, \\ \mathbf{0}, & l \neq 0. \end{cases} \tag{16}$$

*and*

$$\mathbb{E}\left[\boldsymbol{\Psi}_l(\mathbf{X})\boldsymbol{\Psi}_r^T(\mathbf{X})\right] = \begin{cases} \mathbf{I}_{K_{N,l}}, & l = r, \\ \mathbf{0}, & l \neq r. \end{cases} \tag{17}$$

*respectively, where $\mathbf{I}_{K_{N,l}}$ is a $K_{N,l} \times K_{N,l}$ identity matrix.*

*Proof.* Apply the expectation operator on (13) to write $\mathbb{E}[\mathbf{\Psi}_l(\mathbf{X})] = \mathbf{W}_l\mathbb{E}[\mathbf{P}_l(\mathbf{X})]$ and then use (14), with $\mathbf{W}_0 = [1]$ in mind, to derive both the non-trivial and trivial results of (16).

Using (12) in the whitening transformation (13),

$$\begin{aligned}\mathbb{E}[\mathbf{\Psi}_l(\mathbf{X})\mathbf{\Psi}_l(\mathbf{X})^T] &= \mathbf{W}_l\mathbb{E}[\mathbf{P}_l(\mathbf{X})\mathbf{P}_l(\mathbf{X})^T]\mathbf{W}_l^T \\ &= \mathbf{W}_l\mathbf{G}_l\mathbf{W}_l^T \\ &= \mathbf{W}_l\mathbf{W}_l^{-1}\mathbf{W}_l^{-T}\mathbf{W}_l^T = \mathbf{I}_{K_{N,l}},\end{aligned}$$

obtaining the non-trivial result of (17). When $l \neq r$, the trivial result of (17) follows from

$$\mathbb{E}[\mathbf{\Psi}_l(\mathbf{X})\mathbf{\Psi}_r(\mathbf{X})^T] = \mathbf{W}_l\mathbb{E}[\mathbf{P}_l(\mathbf{X})\mathbf{P}_r(\mathbf{X})^T]\mathbf{W}_r^T = \mathbf{0},$$

where the equality to *zero* results from the vanishing expectation as per (15). □

Given the vector $\mathbf{\Psi}_l(\mathbf{x})$ of multivariate orthonormal polynomials, let $\{\Psi_{\mathbf{j}}(\mathbf{x}) : |\mathbf{j}| = l, \mathbf{j} \in \mathbb{N}_0^N\}$, $\Psi_{\mathbf{0}}(\mathbf{x}) = 1$, denote the corresponding set of multivariate orthonormal polynomials. Then $\{\Psi_{\mathbf{j}}(\mathbf{x}) : \mathbf{j} \in \mathbb{N}_0^N\}$ represents an infinite set of multivariate orthonormal polynomials. The second-moment properties follow readily.

**Corollary 10.** *Let $\{\Psi_{\mathbf{j}}(\mathbf{x}) : \mathbf{j} \in \mathbb{N}_0^N\}$ denote an infinite set of multivariate orthonormal polynomials consistent with the probability measure $f_{\mathbf{X}}(\mathbf{x})d\mathbf{x}$. For $\mathbf{j}, \mathbf{k} \in \mathbb{N}_0^N$, the first- and second-order moments of multivariate orthonormal polynomials are*

$$\mathbb{E}\left[\Psi_{\mathbf{j}}(\mathbf{X})\right] = \begin{cases} 1, & \mathbf{j} = \mathbf{0}, \\ 0, & \mathbf{j} \neq \mathbf{0}, \end{cases}$$

*and*

$$\mathbb{E}\left[\Psi_{\mathbf{j}}(\mathbf{X})\Psi_{\mathbf{k}}(\mathbf{X})\right] = \begin{cases} 1, & \mathbf{j} = \mathbf{k}, \\ 0, & \mathbf{j} \neq \mathbf{k}, \end{cases}$$

*respectively.*

## 5. Generalized Polynomial Chaos Expansion

Let $y(\mathbf{X}) := y(X_1, \ldots, X_N)$ be a real-valued, square-integrable output random variable defined on the same probability space $(\Omega, \mathcal{F}, \mathbb{P})$. The vector space $L^2(\Omega, \mathcal{F}, \mathbb{P})$ is a Hilbert space such that

$$\mathbb{E}\left[y^2(\mathbf{X})\right] := \int_{\Omega} y^2(\mathbf{X}(\omega))d\mathbb{P}(\omega) = \int_{\mathbb{A}^N} y^2(\mathbf{x})f_{\mathbf{X}}(\mathbf{x})d\mathbf{x} < \infty$$

with inner product

$$(y(\mathbf{X}), z(\mathbf{X}))_{L^2(\Omega,\mathcal{F},\mathbb{P})} := \int_{\Omega} y(\mathbf{X}(\omega))z(\mathbf{X}(\omega))d\mathbb{P}(\omega) = \int_{\mathbb{A}^N} y(\mathbf{x})z(\mathbf{x})f_{\mathbf{X}}(\mathbf{x})d\mathbf{x} =: (y(\mathbf{x}), z(\mathbf{x}))_{f_{\mathbf{X}}d\mathbf{x}}$$

and norm

$$\|y(\mathbf{X})\|_{L^2(\Omega,\mathcal{F},\mathbb{P})} := \sqrt{(y(\mathbf{X}), y(\mathbf{X}))_{L^2(\Omega,\mathcal{F},\mathbb{P})}} = \sqrt{(y(\mathbf{x}), y(\mathbf{x}))_{f_{\mathbf{X}}d\mathbf{x}}} =: \|y(\mathbf{x})\|_{f_{\mathbf{X}}d\mathbf{x}}.$$

It is elementary to show that $y(\mathbf{X}(\omega)) \in L^2(\Omega, \mathcal{F}, \mathbb{P})$ if and only if $y(\mathbf{x}) \in L^2(\mathbb{A}^N, \mathcal{B}^N, f_{\mathbf{X}}d\mathbf{x})$.

### *5.1. Generalized PCE*

A generalized PCE of a square-integrable random variable $y(\mathbf{X})$ is simply the expansion of $y(\mathbf{X})$ with respect to an orthonormal polynomial basis of $L^2(\Omega, \mathcal{F}, \mathbb{P})$, formally presented as follows.

**Theorem 11.** *Let $\mathbf{X} := (X_1, \ldots, X_N)^T$ be a vector of $N \in \mathbb{N}$ input random variables fulfilling Assumption 1. For $l \in \mathbb{N}_0$, recall that $\boldsymbol{\Psi}_l(\mathbf{x}) \in \mathbb{R}^{K_{N,l}}$, a vector comprising multivariate orthonormal polynomials of degree $l$, is consistent with the probability measure $f_{\mathbf{X}}d\mathbf{x}$. Then*

(1) *for any random variable $y(\mathbf{X}) \in L^2(\Omega, \mathcal{F}, \mathbb{P})$ there exists a Fourier series in multivariate orthonormal polynomials in $\mathbf{X}$, referred to as the generalized PCE of*

$$y(\mathbf{X}) = \sum_{l \in \mathbb{N}_0} \mathbf{C}_l^T \boldsymbol{\Psi}_l(\mathbf{X}), \tag{18}$$

*where the vector of expansion coefficients $\mathbf{C}_l \in \mathbb{R}^{K_{N,l}}$ is defined by*

$$\mathbf{C}_l := \mathbb{E}\left[y(\mathbf{X})\boldsymbol{\Psi}_l(\mathbf{X})\right] := \int_{\mathbb{A}^N} y(\mathbf{x})\boldsymbol{\Psi}_l(\mathbf{x}) f_{\mathbf{X}}(\mathbf{x}) d\mathbf{x}, \;\; l \in \mathbb{N}_0; \tag{19}$$

*and*

(2) *the generalized PCE of $y(\mathbf{X}) \in L^2(\Omega, \mathcal{F}, \mathbb{P})$ converges to $y(\mathbf{X})$ in mean-square, that is, for $y_m(\mathbf{X}) := \sum_{l=0}^{m} \mathbf{C}_l^T \boldsymbol{\Psi}_l(\mathbf{X})$, $m \in \mathbb{N}_0$,*

$$\lim_{m\to\infty} \mathbb{E}\left[|y_m(\mathbf{X}) - y(\mathbf{X})|^2\right] = 0;$$

*converges in probability, that is, for any $\epsilon > 0$,*

$$\lim_{m\to\infty} \mathbb{P}\left(|y_m(\mathbf{X}) - y(\mathbf{X})| > \epsilon\right) = 0;$$

*and converges in distribution, that is, for any $\xi \in \mathbb{R}$,*

$$\lim_{m\to\infty} F_m(\xi) = F(\xi)$$

*such that $F_m(\xi) := \mathbb{P}(y_m(\mathbf{X}) \le \xi)$ and $F(\xi) := \mathbb{P}(y(\mathbf{X}) \le \xi)$ are continuous distribution functions.*

*Proof.* If $y(\mathbf{x}) \in L^2(\mathbb{A}^N, \mathcal{B}^N, f_{\mathbf{X}}d\mathbf{x})$, then by Proposition 3, the expansion

$$y(\mathbf{x}) = \sum_{l \in \mathbb{N}_0} \mathrm{proj}_l y(\mathbf{x}), \tag{20}$$

with $\mathrm{proj}_l y(\mathbf{x}) : L^2(\mathbb{A}^N, \mathcal{B}^N, f_{\mathbf{X}}d\mathbf{x}) \to \mathcal{V}_l^N$ denoting the projection operator, can be formed. Since orthonormalization is a linear transformation, with Proposition 3 in mind, $\mathcal{V}_l^N$ is also spanned by $\{\Psi_{\mathbf{j}}(\mathbf{x}) : |\mathbf{j}| = l, \mathbf{j} \in \mathbb{N}_0^N\}$. Consequently,

$$\mathrm{proj}_l y(\mathbf{x}) = \mathbf{C}_l^T \boldsymbol{\Psi}_l(\mathbf{x}). \tag{21}$$

By definition of the random vector $\mathbf{X}$, the sequence $\{\Psi_{\mathbf{j}}(\mathbf{X})\}_{\mathbf{j}\in\mathbb{N}_0^N}$ or $\{\boldsymbol{\Psi}_l(\mathbf{X})\}_{l\in\mathbb{N}_0}$ is a basis of $L^2(\Omega, \mathcal{F}, \mathbb{P})$, inheriting the properties of the basis $\{\Psi_{\mathbf{j}}(\mathbf{x})\}_{\mathbf{j}\in\mathbb{N}_0^N}$ or $\{\boldsymbol{\Psi}_l(\mathbf{x})\}_{l\in\mathbb{N}_0}$ of $L^2(\mathbb{A}^N, \mathcal{B}^N, f_{\mathbf{X}}d\mathbf{x})$. Therefore, (20) and (21) lead to the expansion in (18).

In reference to Proposition 4, recognize that the set of polynomials from the orthogonal sum

$$\bigoplus_{l\in\mathbb{N}_0} \mathrm{span}\{\Psi_{\mathbf{j}}(\mathbf{x}) : |\mathbf{j}| = l, \mathbf{j} \in \mathbb{N}_0^N\} = \{\Psi_{\mathbf{j}}(\mathbf{x}) : \mathbf{j} \in \mathbb{N}_0^N\} = \Pi^N \tag{22}$$

is also dense in $L^2(\mathbb{A}^N, \mathcal{B}^N, f_{\mathbf{X}}d\mathbf{x})$. Therefore, one has the Bessel's inequality [19]

$$\mathbb{E}\left[\sum_{l\in\mathbb{N}_0} \mathbf{C}_l^T \boldsymbol{\Psi}_l(\mathbf{X})\right]^2 \le \mathbb{E}\left[y^2(\mathbf{X})\right],$$

proving that the generalized PCE converges in mean-square or $L^2$. To determine the limit of convergence, invoke again Proposition 4, which implies that the set $\{\Psi_{\mathbf{j}}(\mathbf{x}) : \mathbf{j} \in \mathbb{N}_0^N\}$ is complete in $L^2(\mathbb{A}^N, \mathcal{B}^N, f_{\mathbf{X}} d\mathbf{x})$. Therefore, Bessel's inequality becomes an equality

$$\mathbb{E}\left[\sum_{l\in\mathbb{N}_0} \mathbf{C}_l^T \mathbf{\Psi}_l(\mathbf{X})\right]^2 = \mathbb{E}\left[y^2(\mathbf{X})\right],$$

known as the Parseval identity [19] for a multivariate orthogonal system, for every random variable $y(\mathbf{X}) \in L^2(\Omega, \mathcal{F}, \mathbb{P})$. Furthermore, as the PCE converges in mean-square, it does so in probability. Moreover, as the expansion converges in probability, it also converges in distribution.

Finally, to find the expansion coefficients, define a second moment

$$e_{\text{PCE}} := \mathbb{E}\left[y(\mathbf{X}) - \sum_{l\in\mathbb{N}_0} \mathbf{C}_l^T \mathbf{\Psi}_l(\mathbf{X})\right]^2 \tag{23}$$

of the difference between $y(\mathbf{X})$ and its full PCE. Differentiate both sides of (23) with respect to $\mathbf{C}_l$, $l \in \mathbb{N}_0$, to write

$$\begin{aligned}
\frac{\partial e_{\text{PCE}}}{\partial \mathbf{C}_l} &= \frac{\partial}{\partial \mathbf{C}_l}\mathbb{E}\left[y(\mathbf{X}) - \sum_{r\in\mathbb{N}_0} \mathbf{C}_r^T \mathbf{\Psi}_r(\mathbf{X})\right]^2 \\
&= \mathbb{E}\left[\frac{\partial}{\partial \mathbf{C}_l}\left\{y(\mathbf{X}) - \sum_{r\in\mathbb{N}_0} \mathbf{C}_r^T \mathbf{\Psi}_r(\mathbf{X})\right\}^2\right] \\
&= 2\mathbb{E}\left[\left\{\sum_{r\in\mathbb{N}_0} \mathbf{C}_r^T \mathbf{\Psi}_r(\mathbf{X}) - y(\mathbf{X})\right\}\mathbf{\Psi}_l^T(\mathbf{X})\right] \\
&= 2\left\{\sum_{r\in\mathbb{N}_0} \mathbf{C}_r^T \mathbb{E}\left[\mathbf{\Psi}_l(\mathbf{X})\mathbf{\Psi}_r^T(\mathbf{X})\right] - \mathbb{E}\left[y(\mathbf{X})\mathbf{\Psi}_l^T(\mathbf{X})\right]\right\} \\
&= 2\left\{\mathbf{C}_l^T - \mathbb{E}\left[y(\mathbf{X})\mathbf{\Psi}_l^T(\mathbf{X})\right]\right\}.
\end{aligned} \tag{24}$$

Here, the second, third, fourth, and last lines are obtained by interchanging the differential and expectation operators, performing the differentiation, swapping the expectation and summation operators, and applying Proposition 9, respectively. The interchanges are permissible as the infinite sum is convergent as demonstrated in the preceding paragraph. Setting $\partial e_{\text{PCE}}/\partial \mathbf{C}_l = 0$ in (24) yields (19), completing the proof. □

The expression of the expansion coefficients can also be derived by simply replacing $y(\mathbf{X})$ in (19) with the full PCE and then using Proposition 9. In contrast, the proof given here demonstrates that the PCE coefficients are determined optimally.

Alternatively, the PCE proposed can be expressed in terms of measure-consistent orthogonal polynomials directly, for instance,

$$y(\mathbf{X}) = \sum_{l\in\mathbb{N}_0} \bar{\mathbf{C}}_l^T \mathbf{P}_l(\mathbf{X}) \tag{25}$$

involving new Fourier coefficients

$$\bar{\mathbf{C}}_l := \mathbf{G}_l^{-1}\mathbb{E}\left[y(\mathbf{X})\mathbf{P}_l(\mathbf{X})\right], \ l \in \mathbb{N}_0. \tag{26}$$

The new coefficients are related to the old coefficients by

$$\bar{\mathbf{C}}_l = \mathbf{W}_l^T \mathbf{C}_l.$$

It is elementary to show that the expansions described by (25) and (26) and (18) and (19) are the same. However, the whitening transformation, which yields measure-consistent orthonormal polynomials, facilitates a relatively simpler form of PCE, as expressed by (18) and (19). This also results in concise expressions of the second-moment properties of PCE, to be discussed in Section 5.2.1. Otherwise, there is no reason to favor one expansion over the other.

The generalized PCE in (18) and (19) should not be confused with that of Xiu and Karniadakis [3]. The PCE presented here further extends the applicability of the existing PCE for arbitrary dependent probability distributions of random input. In contrast, the existing PCE, whether classical [1, 2] or generalized [3], still requires independence of random input.

**Corollary 12.** *Given the preamble of Theorem 11, the generalized PCE can also be expressed in terms of multivariate orthonormal polynomials* $\Psi_{\mathbf{j}}(\mathbf{X})$, $\mathbf{j} \in \mathbb{N}_0^N$, *by*

$$y(\mathbf{X}) = \sum_{\mathbf{j} \in \mathbb{N}_0^N} C_{\mathbf{j}} \Psi_{\mathbf{j}}(\mathbf{X})$$

*with the expansion coefficients*

$$C_{\mathbf{j}} = \mathbb{E}\left[y(\mathbf{X})\Psi_{\mathbf{j}}(\mathbf{X})\right], \ \mathbf{j} \in \mathbb{N}_0^N.$$

*Proof.* Follow the proof of Theorem 11 and apply Corollary 10 to obtain the stated result. □

**Corollary 13.** *Let* $\mathbf{X} = (X_1, \ldots, X_N)^T$ *be a vector of independent, but not necessarily identical, input random variables satisfying Assumption 1, with respective marginal density functions* $f_{X_i}(x_i)$, $i = 1, \ldots, N$. *Denote by* $\Psi_{j_i}(x_i)$ *the* $j_i$*th-degree univariate orthonormal polynomial in* $x_i$, *which is obtained consistent with the probability measure* $f_{X_i}(x_i)dx_i$. *Then the proposed generalized PCE reduces to the traditional PCE, yielding*

$$y(\mathbf{X}) = \sum_{\mathbf{j} \in \mathbb{N}_0^N} C_{\mathbf{j}} \prod_{i=1}^{N} \Psi_{j_i}(X_i)$$

*with the expansion coefficients*

$$C_{\mathbf{j}} = \mathbb{E}\left[ y(\mathbf{X}) \prod_{i=1}^{N} \Psi_{j_i}(X_i) \right].$$

Note that the infinite series in (18) does not necessarily converge almost surely to $y(\mathbf{X})$. Furthermore, there is no guarantee that the moments of PCE of order larger than two will converge. These known fundamental limitations of existing PCE persist in the generalized PCE proposed.

### *5.2. Truncation*

The generalized PCE contains an infinite number of orthonormal polynomials or coefficients. In practice, the number must be finite, meaning that the PCE must be truncated. However, there are multiple ways to perform a truncation. A straightforward approach adopted in this work entails retaining polynomial expansion orders less than or equal to $m \in \mathbb{N}_0$. The result is an $m$th-order generalized PCE approximation[3]

$$y_m(\mathbf{X}) = \sum_{l=0}^{m} \mathbf{C}_l^T \boldsymbol{\Psi}_l(\mathbf{X}) \tag{27}$$

[3]The nouns *degree* and *order* associated with the generalized PCE or orthogonal polynomials are used synonymously in the paper.

of $y(\mathbf{X})$, which contains

$$L_{N,m} = \binom{N+m}{m} = \frac{(N+m)!}{N!m!}$$

number of expansion coefficients defined by (19).

It is natural to ask about the approximation quality of (27). Since the set $\{\boldsymbol{\Psi}_l(\mathbf{x}) : l \in \mathbb{N}_0\}$ or $\{\boldsymbol{\Psi}_l(\mathbf{X}) : l \in \mathbb{N}_0\}$ is complete in $L^2(\mathbb{A}^N, \mathcal{B}^N, f_{\mathbf{X}}d\mathbf{x})$ or $L^2(\Omega, \mathcal{F}, \mathbb{P})$, the truncation error $y(\mathbf{X}) - y_m(\mathbf{X})$ is orthogonal to any element of the subspace from which $y_m(\mathbf{X})$ is chosen, as demonstrated below.

**Proposition 14.** *The truncation error $y(\mathbf{X}) - y_m(\mathbf{X})$ is orthogonal to*

$$\bigoplus_{l=0}^{m} span\{\Psi_{\mathbf{j}}(\mathbf{X}) : |\mathbf{j}| = l, \mathbf{j} \in \mathbb{N}_0^N\}. \tag{28}$$

*Moreover, $\mathbb{E}[y(\mathbf{X}) - y_m(\mathbf{X})]^2 \to 0$ as $m \to \infty$.*

*Proof.* Let

$$\bar{y}_m(\mathbf{X}) := \sum_{r=0}^{m} \bar{\mathbf{C}}_r^T \boldsymbol{\Psi}_r(\mathbf{X}),$$

with arbitrary constant vectors $\bar{\mathbf{C}}_r$, $r = 0, 1, \ldots, m$, be any element of the subspace of $L^2(\Omega, \mathcal{F}, \mathbb{P})$ described by (28). Then

$$\begin{aligned} & \mathbb{E}\left[\{y(\mathbf{X}) - y_m(\mathbf{X})\}\bar{y}_m(\mathbf{X})\right] \\ = \ & \mathbb{E}\left[\left\{\sum_{l=m+1}^{\infty} \mathbf{C}_l^T \boldsymbol{\Psi}_l(\mathbf{X})\right\}\left\{\sum_{r=0}^{m} \bar{\mathbf{C}}_r^T \boldsymbol{\Psi}_r(\mathbf{X})\right\}\right] \\ = \ & \sum_{l=m+1}^{\infty} \sum_{r=0}^{m} \mathbf{C}_l^T \mathbb{E}\left[\boldsymbol{\Psi}_l(\mathbf{X})\boldsymbol{\Psi}_r^T(\mathbf{X})\right] \bar{\mathbf{C}}_r \\ = \ & 0, \end{aligned}$$

where the last line follows from Proposition 9, proving the first part of the proposition. For the latter part, the Pythagoras theorem yields

$$\mathbb{E}[\{y(\mathbf{X}) - y_m(\mathbf{X})\}^2] + \mathbb{E}[y_m^2(\mathbf{X})] = \mathbb{E}[y(\mathbf{X})^2].$$

From Theorem 11, $\mathbb{E}[y_m^2(\mathbf{X})] \to \mathbb{E}[y^2(\mathbf{X})]$ as $m \to \infty$. Therefore, $\mathbb{E}[\{y(\mathbf{X}) - y_m(\mathbf{X})\}^2] \to 0$ as $m \to \infty$. ◻

The second part of Proposition 14 entails $L^2$ convergence, which is the same as the mean-square convergence described in Theorem 11. However, an alternative route is chosen for the proof of Proposition 14.

*5.2.1. Output statistics and other probabilistic characteristics*

The $m$th-order generalized PCE approximation $y_m(\mathbf{X})$ can be viewed as a surrogate of $y(\mathbf{X})$. Therefore, relevant probabilistic characteristics of $y(\mathbf{X})$, including its first two moments and probability density function, if it exists, can be estimated from the statistical properties of $y_m(\mathbf{X})$.

Applying the expectation operator on $y_m(\mathbf{X})$ and $y(\mathbf{X})$ in (27) and (18) and imposing Proposition 9, their means

$$\mathbb{E}\left[y_m(\mathbf{X})\right] = \mathbb{E}\left[y(\mathbf{X})\right] = C_0 \tag{29}$$

are the same as the single element of $\mathbf{C}_0 = (C_0)^T = (C_0)$ and independent of $m$. Therefore, the generalized PCE truncated for any value of $m$ yields the exact mean. Nonetheless, $\mathbb{E}[y_m(\mathbf{X})]$ will be referred to as the $m$th-order generalized PCE approximation of the mean of $y(\mathbf{X})$.

Applying the expectation operator again, this time on $[y_m(\mathbf{X}) - C_0]$ and $[y(\mathbf{X}) - C_0]^2$, and employing Proposition 9 results in the variances

$$\mathrm{var}\left[y_m(\mathbf{X})\right] = \sum_{l=1}^{m} \mathbf{C}_l^T \mathbf{C}_l \tag{30}$$

and

$$\mathrm{var}\left[y(\mathbf{X})\right] = \sum_{l\in\mathbb{N}} \mathbf{C}_l^T \mathbf{C}_l$$

of $y_{S,m}(\mathbf{X})$ and $y(\mathbf{X})$, respectively. In (30), the lower limit exceeds the upper limit when $m = 0$, yielding $\mathrm{var}[y_0(\mathbf{X})] = 0$. This is consistent with $y_0(\mathbf{X}) = C_0$ being a constant function producing no variance. Again, $\mathrm{var}[y_m(\mathbf{X})]$ will be referred to as the $m$th-order generalized PCE approximation of the variance of $y(\mathbf{X})$. Clearly, $\mathrm{var}[y_m(\mathbf{X})]$ approaches $\mathrm{var}[y(\mathbf{X})]$, the exact variance of $y(\mathbf{X})$, as $m \to \infty$.

Note that the formulae for the mean and variance in the existing and generalized PCEs are the same, although the respective expansion coefficients involved are not. The reason for the similarity between the formulae is the orthonormalization as per (13). If the polynomials are chosen not to be normalized, then the formula for the variance will contain additional terms, as recently demonstrated by Rahman [9], although for Gaussian input variables.

Being convergent in probability and distribution, the probability density function of $y(\mathbf{X})$, if it exists, can also be estimated by that of $y_m(\mathbf{X})$. However, deriving analytical formula for the density function is difficult in general. In that case, the density can be estimated by Monte Carlo simulation (MCS) of $y_m(\mathbf{X})$. Such simulation should not be confused with crude MCS of $y(\mathbf{X})$, commonly used for producing benchmark results whenever possible. The crude MCS can be expensive or even prohibitive, particularly when the sample size needs to be very large for estimating tail probabilistic characteristics. In contrast, the MCS embedded in the generalized PCE approximation requires evaluations of simple polynomial functions that describe $y_m$. Therefore, a relatively large sample size can be accommodated in the generalized PCE approximation even when $y$ is expensive to evaluate.

*5.2.2. Computational issues*

There are two important computational issues in the $m$th-order generalized PCE approximation proposed. First, multivariate orthonormal polynomials consistent with the input probability measure must be generated. For Gaussian density on $\mathbb{R}^N$ and select densities on the unit ball $\mathbb{B}^N$ or the simplex $\mathbb{T}^N$, measure-consistent orthogonal polynomials can be generated analytically, to be illustrated in Section 6. In addition, the polynomial moment matrix $\mathbf{G}_l$, required to transform the correlated polynomial vector $\mathbf{P}_l$ into an uncorrelated polynomial vector $\boldsymbol{\Psi}_l$, can be constructed analytically. Therefore, measure-consistent orthonormal polynomials for these probability measures can be produced purely analytically. However, for general probability measures, no such analytical solutions exist; instead, numerical approximations are required. For instance, the Gram-Schmidt process can be employed to generate from monomials a sequence of orthogonal polynomials. However, the process is known to be ill-conditioned. Therefore, more stable methods are needed to compute orthogonal polynomials. Moreover, to construct $\mathbf{G}_l$, deriving an analytical formula for the second-moment properties of orthogonal polynomials for arbitrary non-Gaussian measures is nearly impossible. Having said so, these properties, which represent high-dimensional integrals comprising products of orthogonal polynomials, can be determined by writing them as a sum of expectations of monomials $\{\mathbf{X}^{\mathbf{j}}\}$, $0 \le |\mathbf{j}| \le 2m$, where the moments of $\mathbf{X}$ are calculated either analytically, if possible, or by numerical integration. Note that the numerical integration can be performed with an arbitrary precision even when $N$ is large. This is because no generally expensive output function evaluations are involved.

Second, the calculation of the expansion coefficients requires evaluating the expectations $\mathbb{E}[y(\mathbf{X})\boldsymbol{\Psi}_l(\mathbf{X})]$ for $0 \le l \le m$. These expectations, also various $N$-dimensional integrals on $\mathbb{A}^N$, cannot be determined analytically or exactly if $y$ is a general function. Furthermore, for large $N$, a full numerical integration employing an $N$-dimensional tensor product of a univariate quadrature formula is computationally expensive and likely prohibitive. Therefore, alternative means of estimating these expectations or integrals must be pur-sued. One approach entails exploiting smart combinations of low-dimensional numerical integrations, such as sparse-grid quadrature [20] and dimension-reduction integration [21], to approximate a high-dimensional integral. The other approach consists of efficient sampling methods, such as quasi Monte Carlo simulation [22], importance sampling with Monte Carlo [23], and Markov chain Monte Carlo [24], to name a few. When $y$ is obtained via solution of a differential equation, which is often the case in applications, a frequently used approach is stochastic Galerkin method, where the differential equation is projected onto the same truncated set of orthogonal polynomials analytically, resulting in a system of equations for the coefficients. Nonetheless, more research is needed for robust estimation of the expansion coefficients.

*5.2.3. Implementation*

Algorithm 1 describes a procedure for developing an $m$th-order generalized PCE approximation $y_m(\mathbf{X})$ of a general square-integrable function $y(\mathbf{X})$. It includes calculation of the mean and variance of $y_m(\mathbf{X})$.

---

**Algorithm 1:** Generalized PCE approximation and second-moment statistics

---

**Input**: The total number $N$ of random input variables $\mathbf{X} = (X_1, \ldots, X_N)^T$, a joint probability density function $f_{\mathbf{X}}(\mathbf{x})$ of $\mathbf{X}$ satisfying Assumption 1, a square-integrable function $y(\mathbf{X})$, and the largest order $m$ of orthogonal polynomials

**Output**: The $m$th-order generalized PCE approximation $y_m(\mathbf{X})$ of $y(\mathbf{X})$, and the mean and variance of $y_m(\mathbf{X})$

**for** $l \leftarrow 0$ ***to*** $m$ **do**

  Generate the set of orthogonal polynomials $\{P_{\mathbf{j}}(\mathbf{x}) : |\mathbf{j}| = l, \mathbf{j} \in \mathbb{N}_0^N\}$ consistent with the probability measure $f_{\mathbf{X}}(\mathbf{x})d\mathbf{x}$ of $\mathbf{X}$
  `/* from the Gram-Schmidt process or other means */`

  Construct the orthogonal polynomial vector $\mathbf{P}_l(\mathbf{x})$
  `/* from (5) */`

  Calculate or estimate the polynomial moment matrix $\mathbf{G}_l$
  `/* from (11) analytically or numerically */`

  Perform whitening transformation to produce the orthonormal polynomial vector $\mathbf{\Psi}_l(\mathbf{x})$
  `/* from (13) and Table 1 */`

  Calculate or estimate $\mathbb{E}[y(\mathbf{X})\mathbf{\Psi}_l(\mathbf{X})]$
  `/* from reduced integration or sampling methods */`

  Calculate the vector of expansion coefficients $\mathbf{C}_l$
  `/* from (19) */`

Compile a set $\{\mathbf{C}_l, 0 \leq l \leq m\}$ of at most $m$th-order PCE coefficients and hence construct the $m$th-order PCE approximation $y_m(\mathbf{X})$
`/* from (27) */`

Calculate the mean $\mathbb{E}[y_m(\mathbf{X})]$ and variance $\text{var}[y_m(\mathbf{X})]$
`/* from (29) and (30) */`

---

*5.3. Infinitely many input variables*

In many fields, such as uncertainty quantification, information theory, and stochastic process, functions depending on a countable sequence $\{X_i\}_{i\in\mathbb{N}}$ of input random variables need to be considered. In this case, does the generalized PCE still apply as in the case of finitely many random variables? The answer is yes under certain assumptions, as demonstrated by the following proposition.

**Proposition 15.** *Let $\{X_i\}_{i\in\mathbb{N}}$ be a countable sequence of input random variables defined on the probability space $(\Omega, \mathcal{F}_\infty, \mathbb{P})$, where $\mathcal{F}_\infty := \sigma(\{X_i\}_{i\in\mathbb{N}})$ is the associated $\sigma$-algebra generated. If the sequence $\{X_i\}_{i\in\mathbb{N}}$ satisfies Assumption 1, then the generalized PCE of $y(\{X_i\}_{i\in\mathbb{N}}) \in L^2(\Omega, \mathcal{F}_\infty, \mathbb{P})$, where $y : \mathbb{A}^{\mathbb{N}} \to \mathbb{R}$, converges to $y(\{X_i\}_{i\in\mathbb{N}})$ in mean-square. Moreover, the generalized PCE converges in probability and in distribution.*

*Proof.* According to Proposition 4, $\Pi^N$ is dense in $L^2(\mathbb{A}^N, \mathcal{B}^N, f_{\mathbf{X}} d\mathbf{x})$ and hence in $L^2(\Omega, \mathcal{F}_N, \mathbb{P})$ for every $N \in \mathbb{N}$, where $\mathcal{F}_N := \sigma(\{X_i\}_{i=1}^N)$ is the associated $\sigma$-algebra generated by $\{X_i\}_{i=1}^N$. Here, with a certain abuse of notation, $\Pi^N$ is used as a set of polynomial functions of both real variables $\mathbf{x}$ and random variables $\mathbf{X}$. Now, apply Theorem 3.8 of Ernst et al. [12], which says that if $\Pi^N$ is dense in $L^2(\Omega, \mathcal{F}_N, \mathbb{P})$ for every $N \in \mathbb{N}$, then

$$\Pi^\infty := \bigcup_{N=1}^{\infty} \Pi^N,$$

a subspace of $L^2(\Omega, \mathcal{F}_\infty, \mathbb{P})$, is also dense in $L^2(\Omega, \mathcal{F}_\infty, \mathbb{P})$. But, using (22),

$$\Pi^\infty = \bigcup_{N=1}^{\infty} \bigoplus_{l \in \mathbb{N}_0} \operatorname{span}\{\Psi_{\mathbf{j}} : |\mathbf{j}| = l, \mathbf{j} \in \mathbb{N}_0^N\} = \bigcup_{N=1}^{\infty} \operatorname{span}\{\Psi_{\mathbf{j}} : \mathbf{j} \in \mathbb{N}_0^N\},$$

demonstrating that the set of polynomials from the union is dense in $L^2(\Omega, \mathcal{F}_\infty, \mathbb{P})$. Therefore, the generalized PCE of $y(\{X_i\}_{i\in\mathbb{N}}) \in L^2(\Omega, \mathcal{F}_\infty, \mathbb{P})$ converges to $y(\{X_i\}_{i\in\mathbb{N}})$ in mean-square. Since the mean-square convergence is stronger than the convergence in probability or in distribution, the latter modes of convergence follow readily. □

## 6. Example

### *6.1. Stochastic ODE*

Consider a stochastic boundary-value problem, described by the ordinary differential equation (ODE)

$$-\frac{d}{d\xi}\left(\exp(X_1)\frac{d}{d\xi}y(\xi;\mathbf{X})\right) = \exp(X_2),\ 0 \le \xi \le 1,\ y(\xi;\mathbf{X}) \in \mathbb{R}, \tag{31}$$

with boundary conditions

$$y(0;\mathbf{X}) = 0,\ \exp(X_1)\frac{dy}{d\xi}(1;\mathbf{X}) = 1,$$

where $\mathbf{X} = (X_1, X_2)^T$ is a real-valued, bivariate input random vector with known probability density function. Originally studied by Ernst et al. [12] for a single random variable, the ODE is slightly modified here by introducing two statistically dependent random variables.

Three distinct cases of the probability density function of $\mathbf{X}$, one with an unbounded support and the other two with bounded supports, were considered: (1) a Gaussian density function on $\mathbb{R}^2 := \{(x_1, x_2) : -\infty < x_1, x_2 < +\infty\}$; (2) a rotationally-invariant density function on the unit disk $\mathbb{B}^2 := \{(x_1, x_2) : x_1^2 + x_2^2 \le 1\}$; and (3) a density function on the triangle $\mathbb{T}^2 := \{(x_1, x_2) : 0 \le x_1, x_2; x_1 + x_2 \le 1\}$. For each case of the density function, four subcases, depending on the values of the respective parameters, were studied. Table 2 presents explicit forms of the density functions, including descriptions of the subcases. The values of parameters were chosen to produce widely varying density functions in all three cases, as depicted in Figures 1(a), 2(a), and 3(a). Therefore, there are 12 different density functions of $\mathbf{X}$ in this problem, each leading to a distinct result. The objective is to assess the approximation quality of the truncated generalized PCE in terms of the second-moment statistics of the solution of the ODE for all 12 density functions.

### *6.2. Orthogonal and orthonormal polynomials*

All three input density functions satisfy Items (1) and (2) of Assumption 1. Therefore, measure-consistent orthogonal polynomial bases exist in all cases. However, there are multiple and explicit forms of orthogonal polynomial bases. Indeed, multivariate Hermite polynomials, originally studied by Hermite, are consistent with the Gaussian density function on $\mathbb{R}^N$ as explained by Erdelyi [13]. Recently, the author derived analytical formulae for their first- and second-moment properties, leading to Gauss-Hermite polynomial expansion [9] of a multivariate function. Appell and de Fériet [15], Erdelyi [13], and Dunkl and Xu [10] describe multivariate orthogonal polynomials consistent with the density functions on the unit ball $\mathbb{B}^N$

Table 2: Input probability density functions and measure-consistent orthogonal polynomials.

| Case | Probability density function $f_{X_1X_2}(x_1,x_2)$ | Orthogonal polynomial $P_{j_1j_2}(x_1,x_2)$ |
|---|---|---|
| 1[(a)] | Gaussian density on $\mathbb{R}^2$ : $\frac{1}{2\pi\sigma_1\sigma_2\sqrt{1-\rho^2}}\times \exp\left[-\frac{(\frac{x_1}{\sigma_1})^2-2\rho(\frac{x_1}{\sigma_1})(\frac{x_2}{\sigma_2})+(\frac{x_2}{\sigma_2})^2}{2(1-\rho^2)}\right]$ $0<\sigma_1,\sigma_2<\infty;\ -1<\rho<+1;$ $-\infty<x_1,x_2<+\infty.$ Subcases : $\sigma_1=\sigma_2=1/4,\ \rho=-9/10;$ $\sigma_1=\sigma_2=1/4,\ \rho=-1/2;$ $\sigma_1=\sigma_2=1/4,\ \rho=+1/2;$ $\sigma_1=\sigma_2=1/4,\ \rho=+9/10.$ | $\dfrac{(-1)^{\lvert j_1+j_2\rvert}\dfrac{\partial^{j_1+j_2}}{\partial x_1^{j_1}\partial x_2^{j_2}}f_{X_1X_2}(x_1,x_2)}{f_{X_1X_2}(x_1,x_2)};$ $0\le j_1,j_2<\infty.$ |
| 2[(b)] | Density on the disk $\mathbb{B}^2$ : $\frac{1}{\pi}(\mu+\frac{1}{2})(1-x_1^2-x_2^2)^{\mu-\frac{1}{2}};$ $\mu>-\frac{1}{2};\ (x_1,x_2):x_1^2+x_2^2\le 1.$ Subcases : $\mu=0;\ \mu=1;\ \mu=2;\ \mu=4.$ | $\dfrac{\dfrac{\partial^{j_1+j_2}}{\partial x_1^{j_1}\partial x_2^{j_2}}(1-x_1^2-x_2^2)^{j_1+j_2+\mu-\frac{1}{2}}}{(1-x_1^2-x_2^2)^{\mu-\frac{1}{2}}};$ $0\le j_1,j_2<\infty.$ |
| 3[(b)] | Density on the triangle $\mathbb{T}^2$ : $\frac{\Gamma(\alpha+\beta+\gamma+3)x_1^{\alpha}x_2^{\beta}(1-x_1-x_2)^{\gamma}}{\Gamma(\alpha+1)\Gamma(\beta+1)\Gamma(\gamma+1)};$ $\alpha,\beta,\gamma>-1;$ $(x_1,x_2):0\le x_1,x_2;x_1+x_2\le 1.$ Subcases : $\alpha=\beta=\gamma=0;\ \alpha=\beta=\gamma=1;$ $\alpha=\beta=\gamma=2;\ \alpha=\beta=\gamma=3.$ | $\dfrac{\dfrac{\partial^{j_1+j_2}}{\partial x_1^{j_1}\partial x_2^{j_2}}x_1^{j_1+\alpha}x_2^{j_2+\beta}(1-x_1-x_2)^{j_1+j_2+\gamma}}{x_1^{\alpha}x_2^{\beta}(1-x_1-x_2)^{\gamma}};$ $0\le j_1,j_2<\infty.$ |

(a) See Rahman [9] and/or Erdelyi [13].

(b) See Dunkl and Xu [10].

and on the simplex $\mathbb{T}^N$. In all three cases, a set of measure-consistent orthogonal polynomials $\{P_{\mathbf{j}}(\mathbf{x});\mathbf{j}\in\mathbb{N}_0^N\}$ can be obtained from a multivariate analog of the respective Rodrigues' formula. The third column of Table 2 lists such formulae for $N=2$, which were employed to generate the orthogonal polynomials $\{P_{j_1j_2}(x_1,x_2);j_1,j_2\in\mathbb{N}_0\}$ in this paper. More explicitly, Table 3 presents zeroth-, first-, second-, and third-order (-degree) orthogonal polynomials obtained for the respective last subcases of all three density functions. It is easy to verify from Proposition 8 that all non-zero-degree polynomials have *zero* means and any two distinct polynomials are orthogonal whenever the degrees are different. However, two polynomials of the same degree obtained in this way may not be mutually orthogonal.

Given the orthogonal polynomials, polynomial moment matrices were constructed analytically, which is possible to do for all cases of input density functions employed. Henceforth, the orthonormal polynomials $\Psi_{j_1j_2}(x_1,x_2)$ were generated, also analytically, using the whitening matrix from the Cholesky decomposition. Table 4 lists zeroth-, first-, second-, and third-order orthonormal polynomials calculated for the respective last subcases of the density functions. According to Proposition 9, all non-zero-degree polynomials also have *zero* means, *unit* variances, and any two distinct polynomials are orthogonal whether or not the degrees are the same. Either polynomial set from Tables 3 and 4 can be used to construct a PCE approximation. However, in this paper, orthonormal polynomials from the latter table were used for building PCE approximations and producing subsequent results.

Table 3: A few orthogonal polynomials consistent with the respective last subcases of all three density functions.

Case 1, Subcase 4: Gaussian density on $\mathbb{R}^2$ ($\sigma_1 = \sigma_2 = 1/4, \rho = 9/10$)

$$P_{0,0} = 1,$$
$$P_{1,0} = \frac{160}{19}(10x_1 - 9x_2),$$
$$P_{0,1} = -\frac{160}{19}(9x_1 - 10x_2),$$
$$P_{2,0} = \frac{1600}{361}\left(1600x_1^2 - 2880x_2x_1 + 1296x_2^2 - 19\right),$$
$$P_{1,1} = -\frac{160}{361}\left(14400x_1^2 - 28960x_2x_1 + 14400x_2^2 - 171\right),$$
$$P_{0,2} = \frac{1600}{361}\left(1296x_1^2 - 2880x_2x_1 + 1600x_2^2 - 19\right),$$
$$P_{3,0} = \frac{256000(10x_1 - 9x_2)\left(1600x_1^2 - 2880x_2x_1 + 1296x_2^2 - 57\right)}{6859},$$
$$P_{2,1} = -\frac{51200\left(72000x_1^3 - 209600x_2x_1^2 + 45\left(4496x_2^2 - 57\right)x_1 - 64800x_2^3 + 2489x_2\right)}{6859},$$
$$P_{1,2} = \frac{51200\left(64800x_1^3 - 202320x_2x_1^2 + 131\left(1600x_2^2 - 19\right)x_1 + 45x_2\left(57 - 1600x_2^2\right)\right)}{6859},$$
$$P_{0,3} = -\frac{256000(9x_1 - 10x_2)\left(1296x_1^2 - 2880x_2x_1 + 1600x_2^2 - 57\right)}{6859}.$$

Case 2, Subcase 4: Density on the disk $\mathbb{B}^2$ ($\mu = 4$)

$$P_{0,0} = 1,$$
$$P_{1,0} = -9x_1,$$
$$P_{0,1} = -9x_2,$$
$$P_{2,0} = 11\left(10x_1^2 + x_2^2 - 1\right),$$
$$P_{1,1} = 99x_1x_2,$$
$$P_{0,2} = 11\left(x_1^2 + 10x_2^2 - 1\right),$$
$$P_{3,0} = -429x_1\left(4x_1^2 + x_2^2 - 1\right),$$
$$P_{2,1} = -143x_2\left(10x_1^2 + x_2^2 - 1\right),$$
$$P_{1,2} = -143x_1\left(x_1^2 + 10x_2^2 - 1\right),$$
$$P_{0,3} = -429x_2\left(x_1^2 + 4x_2^2 - 1\right).$$

Case 3, Subcase 4: Density on the triangle $\mathbb{T}^2$ ($\alpha = \beta = \gamma = 3$)

$$P_{0,0} = 1,$$
$$P_{1,0} = 4(-2x_1 - x_2 + 1),$$
$$P_{0,1} = 4(-x_1 - 2x_2 + 1),$$
$$P_{2,0} = 10\left(9x_1^2 + 9(x_2 - 1)x_1 + 2(x_2 - 1)^2\right),$$
$$P_{1,1} = 4\left(9x_1^2 + (23x_2 - 13)x_1 + 9x_2^2 - 13x_2 + 4\right),$$
$$P_{0,2} = 10\left(2x_1^2 + (9x_2 - 4)x_1 + 9x_2^2 - 9x_2 + 2\right),$$
$$P_{3,0} = 60\left(-22x_1^3 - 33(x_2 - 1)x_1^2 - 15(x_2 - 1)^2x_1 - 2(x_2 - 1)^3\right),$$
$$P_{2,1} = 20\left(-22x_1^3 + (42 - 69x_2)x_1^2 - 3\left(17x_2^2 - 25x_2 + 8\right)x_1 - 2(x_2 - 1)^2(5x_2 - 2)\right),$$
$$P_{1,2} = 20\left(-10x_1^3 + (24 - 51x_2)x_1^2 - 3\left(23x_2^2 - 25x_2 + 6\right)x_1 - 22x_2^3 + 42x_2^2 - 24x_2 + 4\right),$$
$$P_{0,3} = 60\left(-2x_1^3 + (6 - 15x_2)x_1^2 + \left(-33x_2^2 + 30x_2 - 6\right)x_1 - 22x_2^3 + 33x_2^2 - 15x_2 + 2\right).$$

Table 4: A few orthonormal polynomials consistent with the respective last subcases of all three density functions.

Case 1, Subcase 4: Gaussian density on $\mathbb{R}^2$ $(\sigma_1 = \sigma_2 = 1/4, \rho = 9/10)$

$$
\begin{aligned}
\Psi_{0,0} &= 1, \\
\Psi_{1,0} &= 4x_1, \\
\Psi_{0,1} &= \frac{40x_2}{\sqrt{19}} - \frac{36x_1}{\sqrt{19}}, \\
\Psi_{2,0} &= 8\sqrt{2}x_1^2 - \frac{1}{\sqrt{2}}, \\
\Psi_{1,1} &= \frac{160x_1x_2}{\sqrt{19}} - \frac{144x_1^2}{\sqrt{19}}, \\
\Psi_{0,2} &= \frac{648}{19}\sqrt{2}x_1^2 - \frac{1440}{19}\sqrt{2}x_1x_2 + \frac{800}{19}\sqrt{2}x_2^2 - \frac{1}{\sqrt{2}}, \\
\Psi_{3,0} &= 32\sqrt{\frac{2}{3}}x_1^3 - 2\sqrt{6}x_1, \\
\Psi_{2,1} &= -288\sqrt{\frac{2}{19}}x_1^3 + 320\sqrt{\frac{2}{19}}x_2x_1^2 + 18\sqrt{\frac{2}{19}}x_1 - 20\sqrt{\frac{2}{19}}x_2, \\
\Psi_{1,2} &= \frac{2592}{19}\sqrt{2}x_1^3 - \frac{5760}{19}\sqrt{2}x_2x_1^2 + \frac{3200}{19}\sqrt{2}x_2^2x_1 - 2\sqrt{2}x_1, \\
\Psi_{0,3} &= -\frac{7776}{19}\sqrt{\frac{6}{19}}x_1^3 + \frac{25920}{19}\sqrt{\frac{6}{19}}x_2x_1^2 - \frac{28800}{19}\sqrt{\frac{6}{19}}x_2^2x_1 + 18\sqrt{\frac{6}{19}}x_1 + \frac{32000}{19}\sqrt{\frac{2}{57}}x_2^3 - 20\sqrt{\frac{6}{19}}x_2.
\end{aligned}
$$

Case 2, Subcase 4: Density on the disk $\mathbb{B}^2$ $(\mu = 4)$

$$
\begin{aligned}
\Psi_{0,0} &= 1, \\
\Psi_{1,0} &= -\sqrt{11}x_1, \\
\Psi_{0,1} &= -\sqrt{11}x_2, \\
\Psi_{2,0} &= \frac{11}{2}\sqrt{\frac{13}{5}}x_1^2 - \frac{1}{2}\sqrt{\frac{13}{5}}, \\
\Psi_{1,1} &= \sqrt{143}x_1x_2, \\
\Psi_{0,2} &= \frac{1}{6}\sqrt{\frac{143}{5}}x_1^2 + \frac{1}{3}\sqrt{715}x_2^2 - \frac{1}{6}\sqrt{\frac{143}{5}}, \\
\Psi_{3,0} &= \frac{3\sqrt{11}x_1}{2} - \frac{13}{2}\sqrt{11}x_1^3, \\
\Psi_{2,1} &= \frac{1}{2}\sqrt{\frac{55}{2}}x_2 - \frac{13}{2}\sqrt{\frac{55}{2}}x_1^2x_2, \\
\Psi_{1,2} &= -\frac{1}{2}\sqrt{\frac{143}{3}}x_1^3 - 5\sqrt{\frac{143}{3}}x_2^2x_1 + \frac{1}{2}\sqrt{\frac{143}{3}}x_1, \\
\Psi_{0,3} &= -\sqrt{\frac{1430}{3}}x_2^3 - \frac{1}{2}\sqrt{\frac{715}{6}}x_1^2x_2 + \frac{1}{2}\sqrt{\frac{715}{6}}x_2.
\end{aligned}
$$

Case 3, Subcase 4: Density on the triangle $\mathbb{T}^2$ $(\alpha = \beta = \gamma = 3)$

$$
\begin{aligned}
\Psi_{0,0} &= 1, \\
\Psi_{1,0} &= \sqrt{\frac{13}{2}} - 3\sqrt{\frac{13}{2}}x_1, \\
\Psi_{0,1} &= -\sqrt{\frac{39}{2}}x_1 - \sqrt{78}x_2 + \sqrt{\frac{39}{2}}, \\
\Psi_{2,0} &= \frac{91x_1^2}{2} - \frac{65x_1}{2} + 5, \\
\Psi_{1,1} &= \frac{21}{2}\sqrt{13}x_1^2 + 21\sqrt{13}x_2x_1 - \frac{27\sqrt{13}x_1}{2} - 6\sqrt{13}x_2 + 3\sqrt{13}, \\
\Psi_{0,2} &= \sqrt{182}x_1^2 + 9\sqrt{\frac{91}{2}}x_2x_1 - 2\sqrt{182}x_1 + 9\sqrt{\frac{91}{2}}x_2^2 - 9\sqrt{\frac{91}{2}}x_2 + \sqrt{182}, \\
\Psi_{3,0} &= -28\sqrt{\frac{221}{3}}x_1^3 + \frac{21}{2}\sqrt{663}x_1^2 - \frac{7\sqrt{663}x_1}{2} + \sqrt{\frac{221}{3}}, \\
\Psi_{2,1} &= -12\sqrt{\frac{4641}{11}}x_1^3 - 24\sqrt{\frac{4641}{11}}x_2x_1^2 + \frac{39}{2}\sqrt{\frac{4641}{11}}x_1^2 + 15\sqrt{\frac{4641}{11}}x_2x_1 - \frac{17}{2}\sqrt{\frac{4641}{11}}x_1 - 2\sqrt{\frac{4641}{11}}x_2 + \sqrt{\frac{4641}{11}}, \\
\Psi_{1,2} &= \quad -4\sqrt{\frac{3094}{3}}x_1^3 - 6\sqrt{9282}x_2x_1^2 + 3\sqrt{9282}x_1^2 - 6\sqrt{9282}x_2^2x_1 + 15\sqrt{\frac{4641}{2}}x_2x_1 - 2\sqrt{9282}x_1 \\
&\quad +3\sqrt{\frac{4641}{2}}x_2^2 - 3\sqrt{\frac{4641}{2}}x_2 + \sqrt{\frac{3094}{3}}, \\
\Psi_{0,3} &= \quad -2\sqrt{\frac{3094}{11}}x_1^3 - 15\sqrt{\frac{3094}{11}}x_2x_1^2 + 6\sqrt{\frac{3094}{11}}x_1^2 - 3\sqrt{34034}x_2^2x_1 + 30\sqrt{\frac{3094}{11}}x_2x_1 - 6\sqrt{\frac{3094}{11}}x_1 \\
&\quad -2\sqrt{34034}x_2^3 + 3\sqrt{34034}x_2^2 - 15\sqrt{\frac{3094}{11}}x_2 + 2\sqrt{\frac{3094}{11}}.
\end{aligned}
$$

### 6.3. Result

#### 6.3.1. Exact solution

A straightforward integration of (31) leads to the exact solution:

$$y(\xi;\mathbf{X}) = \frac{1}{\exp(X_1)}\left[\xi + \left(\xi - \frac{\xi^2}{2}\right)\exp(X_2)\right]. \tag{32}$$

Clearly, the first two raw moments $\mathbb{E}[y(\xi;\mathbf{X})]$ and $\mathbb{E}[y^2(\xi;\mathbf{X})]$, or any probabilistic characteristics of $y(\xi;\mathbf{X})]$ for that matter, depend on the probability density of $\mathbf{X}$. Appendix A provides analytical results of these two moments at $\xi = 1$ for all three cases of input density functions, including those for the individual subcases.

#### 6.3.2. Approximate solution

The Gaussian density function, which has an unbounded support, satisfies Item 3(b) of Assumption 1 [9], whereas the density functions on the unit disk and the triangle, which have bounded supports, clearly fulfill Items 3(b) of Assumption 1. Therefore, the generalized PCE can be applied to solve this problem for all density functions. However, since $y(\xi;\mathbf{X})$ is a non-polynomial function of $\mathbf{X}$, a convergence analysis with respect to $m$ — the order of the generalized PCE approximation — is essential. Employing $m = 1,2,3,4,5,6$ in Algorithm 1, six PCE approximations of $y(\xi;\mathbf{X})$ and their second-moment statistics were constructed or calculated for all three density functions.

Define at $\xi = 1$ an $L^1$ error

$$e_m := \frac{|\mathrm{var}[y(1;\mathbf{X})] - \mathrm{var}[y_m(1;\mathbf{X})]|}{\mathrm{var}[y(1;\mathbf{X})]} \tag{33}$$

in the variance, committed by an $m$th-order generalized PCE approximation $y_m(1;\mathbf{X})$ of $y(1;\mathbf{X})$, where $\mathrm{var}[y(1;\mathbf{X})]$ and $\mathrm{var}[y_m(1;\mathbf{X})]$ are exact and approximate variances, respectively. The exact variance was obtained from the first two raw moments in (A.1) through (A.6), depending on the input probability measure, whereas the approximate variance, given $m$, was calculated following Algorithm 1. The expectations involved in evaluating the expansion coefficients were calculated analytically for the Gaussian density function. However, a mixed analytical-numerical integration was needed and performed with high precision to calculate the expansion coefficients for the other two density functions. Therefore, the variances from the PCE approximations and resultant errors were determined exactly or very accurately.

Figure 1(b) presents four plots describing how the error $e_m$ in (33), calculated for each of the four correlation coefficients of the Gaussian density function, decays with respect to $m$. The attenuation rates for all four correlation coefficients are similar, although the errors for negative correlations are larger than those for positive correlations. The dependency of the error on the sign of the correlation coefficient stems from stronger nonlinearity of the function $y$ with respect to the random input for negative correlations than for positive correlations. Similar plots of error analysis for the density functions on the disk and triangle, each with four distinct subcases, are displayed in Figures 2(b) and Figure 3(b), respectively. Nearly exponential convergence is achieved by the generalized PCE approximations for all three density functions.

While the paper focuses on the theoretical contributions to probability and functional approximation, a brief discussion on the practical significance of the work is warranted. First, the generalized PCE entailing polynomials, orthogonal with respect to the original, non-product-type probability density function, is expected to converge faster than the commonly used tensor product PCE in the transformed variables. This is because the measure transformations – with the exception for dependent Gaussian variables, where such transformations are linear – often lead to highly nonlinear output functions of transformed variables. Second, the generalized PCE proposed is particularly beneficial for non-trivial domains, such as a ball or a simplex, where *a priori* measure transformations are complicated or impractical. Third, stable formulae for computing measure-consistent orthogonal polynomials, at least for the special cases considered in the paper, are highly desirable.

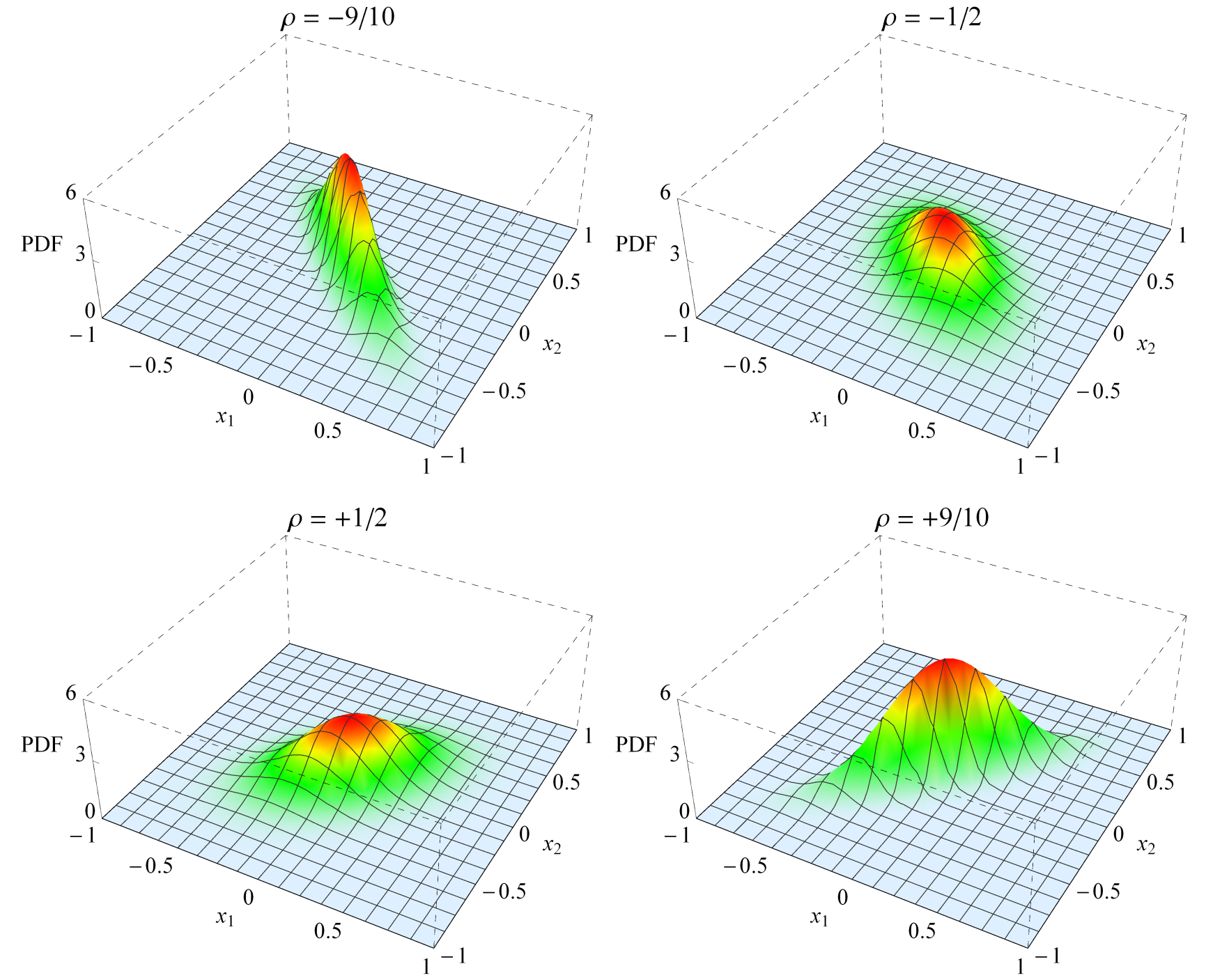

(a)

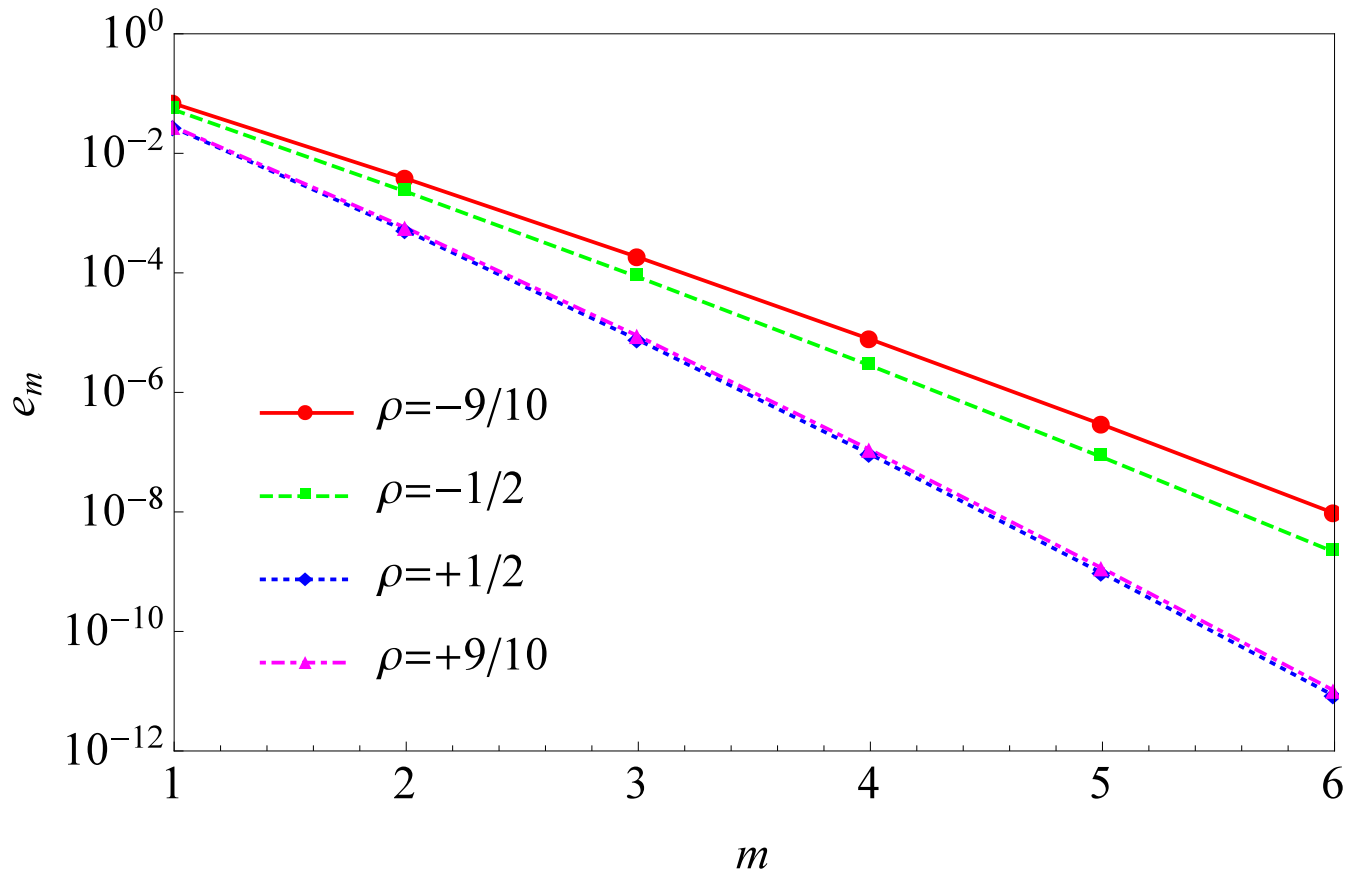

(b)

Figure 1: Input probability measures and PCE results for the Gaussian density on $\mathbb{R}^2$; (a) density of $\mathbf{X}$ for four subcases; (b) decay of $L^1$ error in the variance of $y_m(1;\mathbf{X})$ with respect to $m$.

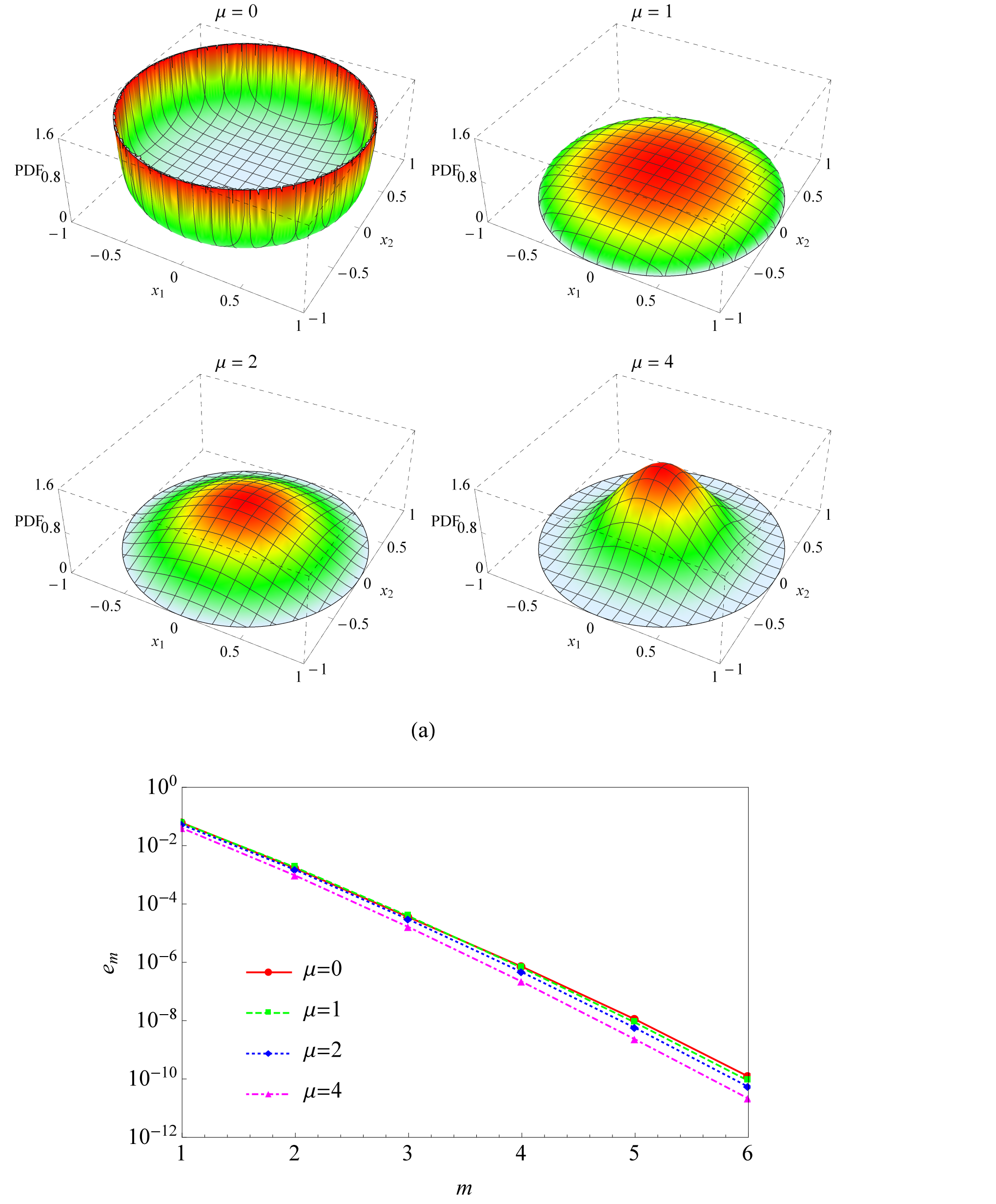

(a)

(b)

Figure 2: Input probability measures and PCE results for the density on disk $\mathbb{B}^2$ ; (a) density of $\mathbf{X}$ for four subcases; (b) decay of $L^1$ error in the variance of $y_m(1;\mathbf{X})$ with respect to $m$.

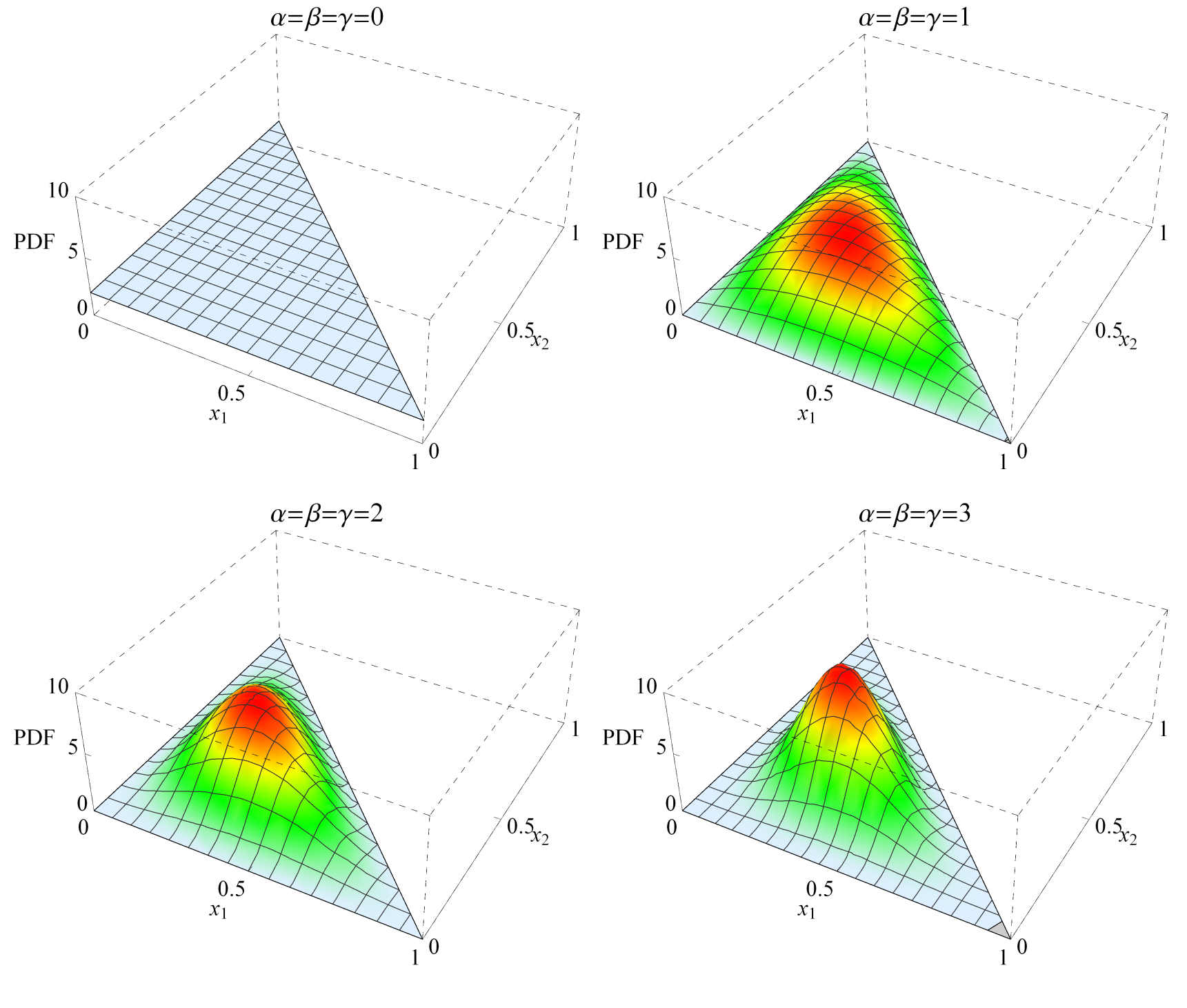

(a)

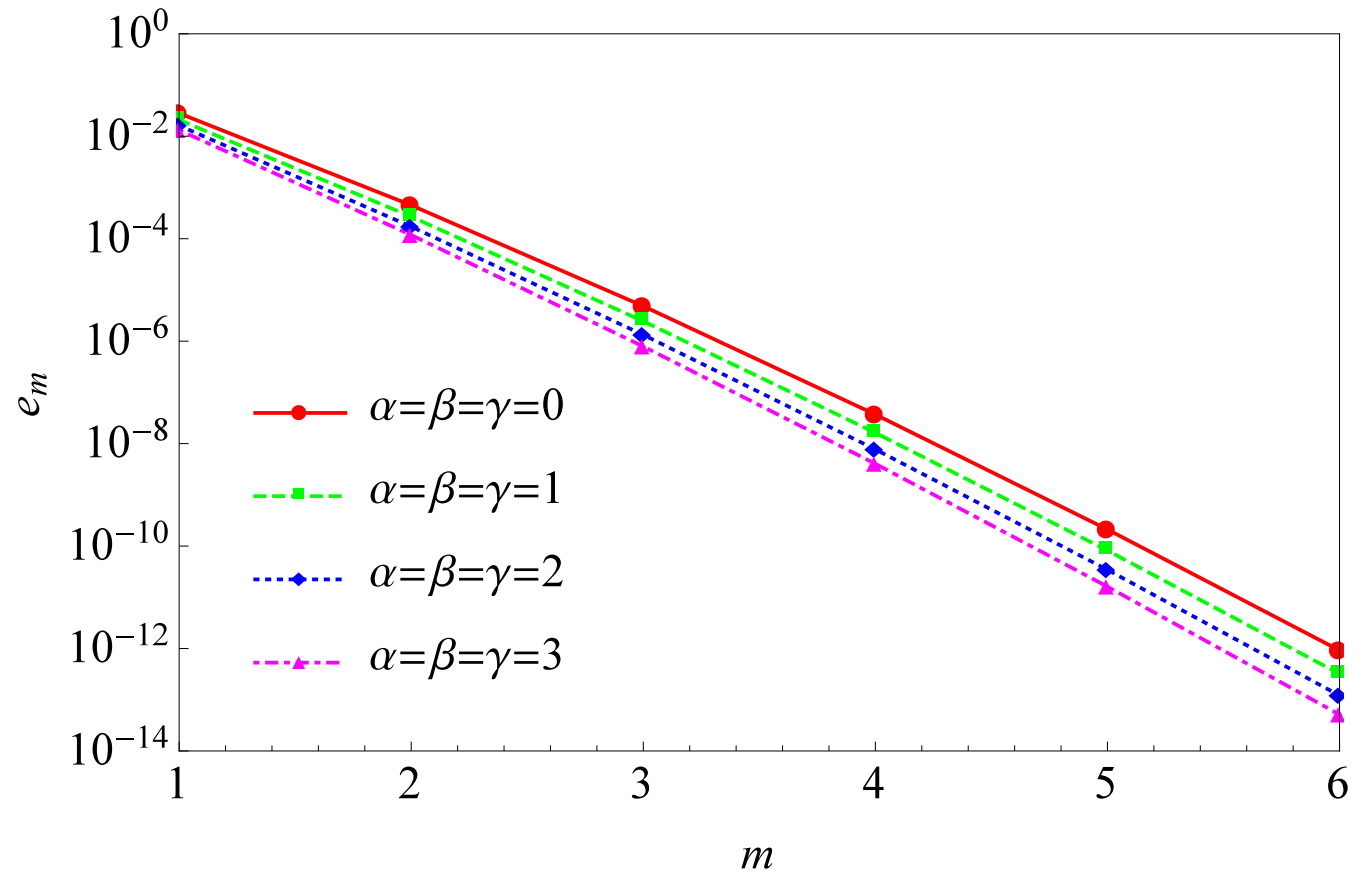

(b)

Figure 3: Input probability measures and PCE results for the density on triangle $\mathbb{T}^2$ ; (a) density of $\mathbf{X}$ for four subcases; (b) decay of $L^1$ error in the variance of $y_m(1;\mathbf{X})$ with respect to $m$.

## 7. Conclusion

A new generalized PCE of a square-integrable random variable, comprising measure-consistent multivariate orthonormal polynomials in dependent random variables with non-product-type probability measures, is presented. There are two main novelties: First, a degree-wise splitting of the polynomial space of all input random variables into orthogonal subspaces, each spanned by measure-consistent multivariate orthogonal polynomials, was constructed, resulting in the PCE developed without the need for a tensor-product structure. Under prescribed assumptions, the set of measure-consistent orthogonal polynomials was proved to form a basis of each subspace, leading to an orthogonal sum of such sets of basis functions to span the space of all polynomials. Second, a whitening transformation is proposed to decorrelate orthogonal polynomials into orthonormal polynomials for an arbitrary probability measure. The transformation is valid whether or not the orthogonal polynomials of the same degree are mutually orthogonal. The orthogonal sum of measure-consistent polynomials, whether orthogonal or orthonormal, is dense in a Hilbert space of square-integrable functions, leading to mean-square convergence of the generalized PCE to the correct limit, including when there are infinitely many random variables. The optimality of the generalized PCE and the approximation quality due to truncation were demonstrated or discussed. For independent probability measures, the proposed generalized PCE reduces to the existing classical or generalized PCE. Analytical formulae are proposed to calculate the mean and variance of a truncated generalized PCE of a general output variable in terms of the expansion coefficients. An example stemming from a stochastic boundary-value problem illustrates the construction and use of a generalized PCE approximation in estimating the statistical properties of an output variable for 12 distinct non-product-type probability measures of input variables.

## Appendix A. Second-moment properties of $y(1;\mathbf{X})$

Applying the expectation operators on (32) and its square, the first two raw moments of $y(1;\mathbf{X})$ are respectively given by (A.1) and (A.2) for the Gaussian density on $\mathbb{R}^2$, by (A.3) and (A.4) for the density on the unit disk $\mathbb{B}^2$, and by (A.5) and (A.6) for the density on the triangle $\mathbb{T}^2$.

(1) Gaussian density on $\mathbb{R}^2$ ($0<\sigma_1,\sigma_2<\infty,\ -1<\rho<+1$):

$$\mathbb{E}\left[y(1;\mathbf{X})\right] = \frac{1}{2}e^{\frac{\sigma_1^2}{2}}\left[e^{\frac{1}{2}\sigma_2(\sigma_2-2\rho\sigma_1)}+2\right] \approx \begin{cases} 1.59479, & \sigma_1=\sigma_2=\frac{1}{4},\ \rho=-\frac{9}{10}, \\ 1.58089, & \sigma_1=\sigma_2=\frac{1}{4},\ \rho=-\frac{1}{2}, \\ 1.54762, & \sigma_1=\sigma_2=\frac{1}{4},\ \rho=+\frac{1}{2}, \\ 1.53488, & \sigma_1=\sigma_2=\frac{1}{4},\ \rho=+\frac{9}{10}. \end{cases} \tag{A.1}$$

$$\mathbb{E}\left[y^2(1;\mathbf{X})\right] = \frac{1}{4}e^{2\sigma_1^2}\left[4e^{\frac{1}{2}\sigma_2(\sigma_2-4\rho\sigma_1)}+e^{2\sigma_2(\sigma_2-2\rho\sigma_1)}+4\right] \approx \begin{cases} 2.84348, & \sigma_1=\sigma_2=\frac{1}{4},\ \rho=-\frac{9}{10}, \\ 2.74142, & \sigma_1=\sigma_2=\frac{1}{4},\ \rho=-\frac{1}{2}, \\ 2.51472, & \sigma_1=\sigma_2=\frac{1}{4},\ \rho=+\frac{1}{2}, \\ 2.43420, & \sigma_1=\sigma_2=\frac{1}{4},\ \rho=+\frac{9}{10}. \end{cases} \tag{A.2}$$

(2) Density on the unit disk $\mathbb{B}^2$ ($\mu=0,1,2,4$):

$$\mathbb{E}\left[y(1;\mathbf{X})\right] = \begin{cases} \frac{1}{4}\left[4\sinh(1)+\sqrt{2}\sinh(\sqrt{2})\right] \approx 1.85935, & \mu=0, \\ \frac{3}{8e}\left[8-\sqrt{2}e\sinh(\sqrt{2})+2e\cosh(\sqrt{2})\right] \approx 1.71105, & \mu=1, \\ \frac{15}{16e}\left[-56+8e^2+5\sqrt{2}e\sinh(\sqrt{2})-6e\cosh(\sqrt{2})\right] \approx 1.64896, & \mu=2, \\ \frac{945}{64e}\left[-8512+1152e^2+199\sqrt{2}e\sinh(\sqrt{2})-250e\cosh(\sqrt{2})\right] \approx 1.59358, & \mu=4. \end{cases} \tag{A.3}$$

$$
\mathbb{E}\left[y^2(1;\mathbf{X})\right] = \begin{cases} \frac{1}{80}\left[40\sinh(2)+5\sqrt{2}\sinh(2\sqrt{2})+16\sqrt{5}\sinh(\sqrt{5})\right] \approx 4.6268, & \mu=0, \\ \frac{3}{3200e^2}\left[600+200e^4-25\sqrt{2}e^2\sinh(2\sqrt{2})- \right. \\ \left. 128\sqrt{5}e^2\sinh(\sqrt{5})+100e^2\cosh(2\sqrt{2})+640e^2\cosh(\sqrt{5})\right] \approx 3.57604, & \mu=1, \\ -\frac{3}{25600}\left[-28000\sinh(2)-1375\sqrt{2}\sinh(2\sqrt{2})- \right. \\ 8192\sqrt{5}\sinh(\sqrt{5})+24000\cosh(2)+ \\ \left. 1500\cosh(2\sqrt{2})+15360\cosh(\sqrt{5})\right] \approx 3.15891, & \mu=2, \\ -\frac{189}{8192000}\left[-24080000\sinh(2)-330625\sqrt{2}\sinh(2\sqrt{2})- \right. \\ 4653056\sqrt{5}\sinh(\sqrt{5})+23200000\cosh(2)+ \\ \left. 462500\cosh(2\sqrt{2})+10158080\cosh(\sqrt{5})\right] \approx 2.80304. & \mu=4. \end{cases} \tag{A.4}
$$

(3) Density on the triangle $\mathbb{T}^2$ ($\alpha=\beta=\gamma=0,1,2,3$):

$$
\mathbb{E}\left[y(1;\mathbf{X})\right] = \begin{cases} \frac{1}{2e}\left(5-2e+e^2\right) \approx 1.27884, & \alpha=\beta=\gamma=0, \\ -\frac{10}{e}\left(-66+16e+3e^2\right) \approx 1.25198, & \alpha=\beta=\gamma=1, \\ \frac{168}{e}\left(6855-2848e+120e^2\right) \approx 1.24129, & \alpha=\beta=\gamma=2, \\ -\frac{5940}{e}\left(-1092315+381956e+7315e^2\right) \approx 1.23556, & \alpha=\beta=\gamma=3. \end{cases} \tag{A.5}
$$

$$
\mathbb{E}\left[y^2(1;\mathbf{X})\right] = \begin{cases} \frac{1}{24e^2}\left(23-15e^2+16e^3+3e^2\sinh(2)\right) \approx 1.77024, & \alpha=\beta=\gamma=0, \\ -\frac{5}{576e^2}\left(-3751-6840e^2+2560e^3+27e^4\right) \approx 1.64355, & \alpha=\beta=\gamma=1, \\ \frac{7}{288}\left(-1075572+378880e-281645\sinh(2)+ \right. \\ \left. 283670\cosh(2)\right) \approx 1.59375, & \alpha=\beta=\gamma=2, \\ -\frac{55}{27648e^2}\left(-2329128235-16257780864e^2+ \right. \\ \left. 6087639040e^3+3393495e^4\right) \approx 1.56717, & \alpha=\beta=\gamma=3. \end{cases} \tag{A.6}
$$

## Acknowledgments

The author thanks the reviewers for providing useful comments.